\documentclass[12pt]{amsart}
\usepackage{amssymb,amsmath}
\def\R {\mathbb{R}}

\def\Dt{\partial_t}
\def\Ree{\operatorname{Re}}
\def\Imm{\operatorname{Im}}

\def\mm{\rm m}
\def\ff{\rm f}

\newtheorem{proposition}{Proposition}[section]
\newtheorem{theorem}[proposition]{Theorem}
\newtheorem{corollary}[proposition]{Corollary}
\newtheorem{lemma}[proposition]{Lemma}
\theoremstyle{definition}

\newtheorem{remark}[proposition]{Remark}

\numberwithin{equation}{section}

\def \au {\rm}
\def \ti {\it}

\def \no#1#2#3 {{\bf #1} (#3), #2.}
\def \eds#1#2#3 {#1, #2, #3.}
\def\Bbb{\mathbb}
\def\eb{\varepsilon}
\def\({\left(}
\def\){\right)}
\def\xx{\text{\bf x}}
\def\yy{\text{\bf y}}
\def\zz{\text{\bf z}}
\def\pp{\text{\bf p}}
\def\ff{\text{\bf f}}
\def\ggg{\text{\bf g}}
\def\qq{\text{\bf q}}
\def\FF{\text{\bf F}}
\def\GG{\text{\bf G}}
\def\vv{\text{\bf v}}
\def\<{\left<}
\def\>{\right>}

\def\uu{\text{\bf u}}

\def\hh{\text{\bf h}}

\def\sspan{\operatorname{span}}
\def\Cal{\mathcal}
\title [Space-Time chaos in Ginzburg-Landau equations]{Analytical proof of space-time chaos in Ginzburg-Landau equations}
\author{D. Turaev and S. Zelik}

\begin{document}
\begin{abstract} We prove that the attractor of the
1D quintic complex Ginzburg-Landau equation with a broken phase symmetry has strictly positive
space-time entropy for an open set of parameter values. The result is obtained by studying chaotic oscillations
in grids of weakly interacting solitons in a class of Ginzburg-Landau type equations. We provide an analytic proof
for the existence of two-soliton configurations with chaotic temporal behavior, and construct solutions which are closed
to a grid of such chaotic soliton pairs, with every pair in the grid well spatially separated from the neighboring ones
for all time. The temporal evolution of the well-separated multi-soliton structures is described
by a weakly coupled lattice dynamical system (LDS) for the coordinates and phases of the solitons.
We develop a version of normal hyperbolicity theory for the weakly coupled LDS's with continuous time
and establish for them the existence of space-time chaotic patterns similar to the Sinai-Bunimovich chaos in
discrete-time LDS's. While the LDS part of the theory may be of independent interest, the main difficulty
addressed in the paper concerns with lifting the space-time chaotic solutions of the LDS back to
the initial PDE. The equations we consider here are space-time autonomous, i.e. we impose no spatial
or temporal modulation which could prevent the individual solitons in the grid from drifting towards each other
and destroying the well-separated grid structure in a finite time. We however manage to show that the set of space-time
chaotic solutions for which the random soliton drift is arrested is large enough, so the corresponding
space-time entropy is strictly positive.
\end{abstract}
\maketitle
\section{Introduction}\label{si}
We demonstrate that if an evolutionary system of partial differential equations (PDE)
in unbounded domain has a solution localized in space and chaotic in time, then
one should expect both temporal and spatial chaotic behavior in the system. Namely,
one may observe a formation of non-trivial spatial patterns that evolve in an irregular
fashion with time, and the corresponding {\em space-time entropy} \cite{15,56} is strictly positive.
In other words, the number of solutions which are essentially different from each other
on a finite space-time window grows exponentially with the window volume.

As a tool for finding the spatially-localized, temporally-chaotic solutions one may try,
as we do it here, to look for special types of both spatially and temporally localized solutions.
Thus, like Shilnikov homoclinic loop and Lorenz butterfly serve as a criterion for
chaos formation in systems of ODE's \cite{Sh1,Sh2,Sh3,Sh4}, the existence of the Shilnikov homoclinic loop
in the dynamical system generated by the PDE on the space of spatially localized solutions implies
the space-time chaos in the extended system that corresponds to uniformly bounded solutions of the same PDE.

We do not prove this principle in full generality here. Instead, we decided to show
how it works for a class of Ginzburg-Landau equations with a broken phase symmetry.
The main motivation for such
approach is that despite a huge amount of numerical and experimental data on different
types of space-time irregular behavior in various systems,
there are very few rigorous mathematical results on this topic
and mathematically relevant models describing these phenomena. Therefore, we made
an effort of providing a free from numerics, completely analytic proof
of the existence of space-time chaos in an important equation of mathematical physics.

The basic mathematical model for the space-time chaotic behavior
is the so-called Sinai-Bunimovich chaos in discrete lattice dynamics, see \cite{3,14,34,35}.
This model consists of a $\Bbb Z^n$-grid of discrete-time
chaotic oscillators coupled by a weak interaction. The single chaotic oscillator of this grid is described, say,
by the Bernoulli scheme $\Cal M^1:=\{0,1\}^{\Bbb Z}$, so the uncoupled system naturally has an
infinite-dimensional hyperbolic set homeomorphic to the multi-dimensional
Bernoulli scheme $\Cal M^{n+1}:=\{0,1\}^{\Bbb Z^{n+1}}=(\Cal M^1)^{\Bbb Z^n}$.
The temporal evolution operator is then conjugate to
the shift in $\Cal M^{n+1}$ along the first coordinate
and the other $n$ coordinate shifts are associated with the spatial translations on the grid.
Due to the structural stability of hyperbolic sets,
the above structure survives under a sufficiently weak coupling. Thus, in this
model, the space-time chaos is described by
the multi-dimensional Bernoulli scheme $\Cal M^{n+1}$.

Importantly, the space-time entropy in the Sinai-Bunimovich model is strictly positive.
We know from the general theory of dissipative systems in unbounded domains
(see e.g. \cite{15,281,54,55,56}) that under some reasonable dissipativity assumptions
this entropy is finite for systems of evolutionary PDE's, therefore the Sinai-Bunimovich model
carries ``enough complexity'' to be able to capture certain basic features of spatio-temporal
chaos in systems of various nature. In particular, it is well established by now (see e.g. \cite{Po,Man9,BaTu5})
that the transition from regular to chaotic space-time behavior often happens via the emergence
of well spatially separated and long living ``turbulent spots''. As the interaction between such spots seems
to be weak, the Sinai-Bunimovich chaos paradigm can be relevant for the analysis of these near-threshold phenomena.

Yet, a direct application of the Sinai-Bunimovich construction to systems with continuous time and space
is not possible, in general. Even the existence of one PDE
which possesses an infinite-dimensional Bernoulli scheme was a long-standing open problem.
The first examples of such PDEs (in the class of reaction-diffusion systems), have been
recently constructed in \cite{32}. The method used in that
paper is based on a strong and explicit spatio-temporal modulation of the equation right-hand sides,
which effectively transforms the systems into a discrete-time lattice dynamical system. The disadvantage is that
very special (and artificial) nonlinear interaction functions emerge in the result, which are far from the usual
nonlinearities arising in physics models.

A different approach to the problem is suggested in \cite{57}, where a theory of weak interaction
of dissipative solitons was developed and, as an application,
a space-time chaotic pattern has been constructed for the perturbed 1D Swift-Hohenberg equation
\begin{equation}\label{0.15}
\Dt u+(\partial_x^2+1)^2u+\beta^2 u+f(u)=\mu h(t,x,\mu),\ \  f(u)=u^3+\kappa u^2,\ \ \mu\ll1.
\end{equation}
Here $\mu h(t,x,\mu)$ is a space-time periodic forcing. Its exact form is quite non-trivial, however
the amplitude $\mu$ can be taken arbitrarily small. The idea is to create a spatially
localized spot of chaotic temporal behavior, and to build then a grid of such spots, well separated in
space. The spots are pinned down to the prescribed locations at the grid points by spatial oscillations
in the forcing $\mu h$. If the spots stay sufficiently far apart, their interaction is small, so a small
amplitude forcing occurs to be sufficient to sustain the grid for all times (the wave length of the forcing
has, however, to grow as the amplitude decreases).

Equation \eqref{0.15} at $\mu=0$, like many other important equations, does have a spatially localized
solution, {\em a soliton}, $u=U(x)$ with exponentially decaying tails.
One may therefore look, at all small $\mu$, for {\em multi-soliton} solutions in the form
$$
u(t,x)=\sum_jU(x-\xi_j(t))+{\rm ``small~corrections"},
$$
where $\xi_j(t)$ is the position of the $j$-th soliton; the well-separation condition reads as
$L:=\inf_{j\ne k}\|\xi_j-\xi_k\|\gg1$. Due to the ``tail'' interaction and the small forcing, the solitons' positions
$\xi_j(t)$ may move slowly, and this motion is described by a lattice dynamical system (LDS), see \cite{57} for details.
The obtained LDS is not in the form one needs for establishing the Sinai-Bunimovich chaos
(a grid of chaotic maps with weak coupling), since the individual
solitons $u=U(x)$ are {\it equilibria} at $\mu=0$ and do not have their own (chaotic) dynamics.
However, as it is shown in \cite{57}, a pair of weakly interacting solitons
in the 1D Swift-Hohenberg equation can be forced to oscillate chaotically in time by an appropriate
choice of the time-periodic perturbation $\mu h(t,x,\mu)$. For a well-separated grid of such soliton pairs,
one obtains a time-periodic LDS, and the period map for this system is the sought discrete lattice of weakly coupled
chaotic maps, i.e the space-time chaos is established.

The scope of \cite{57} is much more general than the Swift-Hohenberg equation: by developing
the center manifold approach proposed in \cite{San}, the paper derives the LDS
that governs the evolution of weakly coupled multi-soliton configurations for a large class of systems of
evolutionary PDE's. It also proposes a method for constructing spatially localized and temporally chaotic
solutions which are obtained as a system of finitely many weakly coupled stationary solitons. Note
that, although spatially localized solutions with non-trivial temporal dynamics have been observed numerically and
experimentally in various physical systems (see e.g. \cite{055,12,46} and references therein), the direct analytic
detection and study of such solutions is obviously a very difficult task. However, when a finite system of
well-separated solitons is considered, the descriprion provided by \cite{57} for the evolution of such object
is often a low-dimensional system of ODE's which can exhibit a chaotic dynamics \cite{46} and
can be studied analytically, so the chaotic temporal behavior of such localized patterns can be rigorously proven.

In the present paper we show how a space-time chaotic lattice can be built out of these chaotic multi-soliton
systems in the case where {\em no spatial nor temporal modulation} is imposed. Two problems immediately
appear in this setting:\\
1. with no external forcing, the LDS which describes the multi-soliton dynamics is an autonomous system with continuous
time, so the Sinai-Bunimovich chaos construction (for which the discreteness of time is very essential) is not
applicable;\\
2. with no spatial modulation, there is no pinning mechanism which would keep the solitons eternally close
to any given spatial grid, therefore the infinite-time validity of the LDS description is no longer guaranteed.

We resolve here both the issues. As an application, we
consider the 1D quintic complex Ginzburg-Landau equation with slightly broken phase symmetry:
\begin{equation}\label{0.nlsgl}
\Dt u=(1+i\beta)\partial_{x}^2 u-(1+i\delta)u+(i+\rho)|u|^2u-(\varepsilon_1+i\varepsilon_2)|u|^4u+\mu,
\end{equation}
where $\beta,\delta,\rho,\varepsilon_{1,2},\mu$ are some real parameters, and $\mu\ll1$.
We mention that, in contrast to the previous on the Swift-Hohenberg equation,
we do not have here any artificial functions, and the only freedom we have is the choice
of the numeric parameters. Note also that the Ginzburg-Landau equation serves as a normal form
near an onset of instability, i.e. it very often appears in applications
as a modulation equation for various more complicated problems. The phase symmetry in the
modulation equation appears as an artefact of closeness to the instability threshold, so if there
is no such symmetry in the original problem, then the effects of small symmetry breaking
also need to be considered, see \cite{311} and references therein.
While we introduce only the simplest symmetry breaking term (``$+\mu$'') in \eqref{0.nlsgl}, the general case
is also covered by the theory (see Section \ref{s4}).

The main result of the paper is the following theorem (Section \ref{nls}).
\begin{theorem} There exists an open set of parameters $(\beta,\delta,\rho,\eb_1,\eb_2,\mu)$ such that equation \eqref{0.nlsgl}
possesses a global attractor $\mathcal A$ (say, in the phase space $L^2_b(\R)$) with strictly positive space-time entropy
$$
h_{s-t}(\mathcal A)>0.
$$
\end{theorem}

Equation \eqref{0.nlsgl} at $\mu=0$ has the additional phase
symmetry $u\to e^{i\phi}u$. Therefore, for each stationary soliton $u=V(x)$ of this equation,
$u=e^{i\phi} V(x)$ is also a stationary soliton. Therefore, the multi-soliton configurations
are given by
$$
u(t,x)=\sum_{i} e^{i\phi_i(t)}V(x-\xi_i(t))+ {\rm ``small~corrections"},
$$
where $\xi_j$ and $\phi_j$ are the coordinate
and phase of the $j$-th soliton. For a soliton pair with the states
$(\xi_1,\phi_1)$ and $(\xi_2,\phi_2)$, the evolution is governed,
to the leading order with respect to the distance $|\xi_2-\xi_1|$, by the following system of ODE's:
\begin{equation}\label{0.2sol}
\begin{cases}
\frac d{d\tau}R=ae^{-\alpha R}\sin(\omega R+\theta_1)\cos(\Phi),
\\
\frac d{d\tau}\Phi=be^{-\alpha R}\cos(\omega R+\theta_2)\sin(\Phi)
-2c\nu\sin(\frac{\Phi}{2})\sin(\Psi),
\\
\frac d{d\tau}\Psi=\frac{b}{2} e^{-\alpha R}\sin(\omega R+\theta_2)\cos(\Phi)
+c\nu\cos(\frac{\Phi}{2})\cos(\Psi)-\Omega,
\end{cases}
\end{equation}
see \cite{Vla,46,57}. Here $\tau$ is a scaled slow time, $R:=(\xi_2-\xi_1)/2$, $\Phi:=\phi_1-\phi_2$,
$\Psi:=(\phi_1+\phi_2)/2$
and $a$, $b$, $\omega$, $\theta_{1,2}$, $c$, $\nu$ and $\Omega$ are parameters whose exact values
depend on the values of the original parameters
of \eqref{0.nlsgl} (see the corresponding expressions, as well as asymptotic expansions
near the exactly solvable nonlinear Schr\"odinger equation, in Sections \ref{s4},\ref{nls}).
While the variables $R,\Phi$ and $\Psi$ can be treated as the ``internal variables'' of the two-soliton pattern,
the variable $p:=(\xi_1+\xi_2)/2$ marks the spatial position of the soliton pair. To the leading order,
it is governed by the equation
\begin{equation}\label{0.mc}
\frac d{d\tau}p=\frac{a}{2}e^{-\alpha R}\cos(\omega R+\theta_1)\sin(\Phi):=g(R,\Phi).
\end{equation}

A numerical study of system \eqref{0.2sol} undertaken in \cite{46} revealed various chaotic regimes
for different parameter values. In order to provide an analytic proof (see Lemma \ref{Lemzero3})
of the chaotic behavior (i.e. the existence of a nontrivial hyperbolic invariant set)
in this system for an open set of parameter values, we have found a
degenerate equilibrium of that system with 3 zero eigenvalues.
The normal form equation for this bifurcation (see \cite{CTS}) leads to the following 3rd order equation:
$$
Y^{\prime\prime\prime} = 1-Y^2 + E Y^{\prime},
$$
where $E$ (a certain combination of the parameters of the original system) can take any real values.
In \cite{Ibanez}, the existence of a Shilnikov homoclinic loop for this system was proven at certain $E$
values, which implies \cite{Sh1,Sh2} chaos for some interval of the $E$ values and, hence, for an open set
of parameter values for system \eqref{0.2sol}.

Chaotic solutions of system \eqref{0.2sol} correspond to a chaotically oscillating soliton pair, which
is a temporally chaotic and spatially localised solution, by construction.
After that, according to the program described above, we build a well spatially separated
lattice of such time-chaotic solitons.
The center manifold reduction theorem proved in \cite{57} ensures that the evolution of this lattice
is governed by a system of infinitely many weakly coupled copies of the ODE's \eqref{0.2sol},\eqref{0.mc}.

Even when every individual ODE-subsystem in the continuous time LDS is hyperbolic, the LDS itself is not hyperbolic
(this is a principal difference with the Sinai-Bunimovich chaos in the discrete-time LDS's where the
countable product of hyperbolic sets for the individual maps is hyperbolic again). Each constituent ODE contributes
a neutral direction corresponding to the time shift, so for the linearized flow of the
continuous time LDS we have infinitely many neutral directions. Therefore, after a weak coupling is switched on,
the dynamics is not preserved (the LDS can hardly be topologically conjugate to the uncoupled one). Still,
the invariant manifold theorem of Section \ref{s2} shows that if, given any orbit of the uncoupled LDS, we consider
the family of all orbits obtained by all possible time-reparametrizations in each of the constituent ODE's,
then this family continues in a unique way as an invariant manifold of the weakly coupled LDS. This fact allows
to show the strict positivity of space-time topological entropy for the countable systems of weakly coupled chaotic
oscillators with continuous time.

In fact, results of Section \ref{s2} cover LDS's of a more general type. The problem we have to deal with is that,
although system \eqref{0.2sol} for the internal variables $(R,\Phi,\Psi)$ of the chaotic soliton does have
a uniformly hyperbolic set,
the full system describing the motion of the chaotic soliton includes equation \eqref{0.mc} for the soliton position $p$,
and is clearly non-hyperbolic (so we have to consider the LDS's built of partially-hyperbolic individual ODE's).
The neutral directions appear because the right-hand sides of \eqref{0.2sol},\eqref{0.mc} are $p$-independent, which is
a mere consequence of the translational symmetry of the PDE under consideration, i.e. their presence
is an inherent property of the soliton-interaction equations in systems without a spatial modulation.

Since the internal variables of the soliton change chaotically with time, the soliton position $p(t)$
performs, essentially, an unbounded random walk (as an integral of a chaotic input, see \eqref{0.mc}).
When the chaotic solitons are well spatially separated, the contribution of the neighboring solitons
to the $p$-equation is small, so we have essentially independent random walks for each of the chaotic
solitons in the lattice. This makes it impossible for us to ensure that the well-separation condition is fulfilled
for all times and all initial multi-soliton configurations. We, in fact, believe that the majority of
these configurations do break up this condition in a finite time, so the corresponding solutions
cannot be completely described by the weak soliton interaction paradigm.

However, the weak soliton interaction theory of \cite{57} is the only tool we have here
for the analysis of dynamics of the multi-soliton patterns. As we are unable to control the
soliton's random walk, we devise a method of keeping track of those configurations
for which the well-separation condition holds eternally (i.e. the LDS description is applicable).
This method allows us to verify (Section \ref{s3}) that the number of such solutions
is large enough to ensure the positivity of the space-time entropy. It is worth to emphasize that,
instead of fighting with the random walks, our method exploits them in a crucial way.

Roughly speaking, assume that the hyperbolic set for \eqref{0.2sol} contains
two periodic orbits $\Gamma_1$ and $\Gamma_2$ and a number of heteroclinics which connect them.
Assume that, according to equation \eqref{0.mc}, the soliton pair moves
to the left if $(R,\Phi,\Psi)$ belongs to $\Gamma_1$ and to the right when it belong to $\Gamma_2$.
The direction of this motion is determined by the sign of
$\displaystyle b_j:=\frac{1}{T_j}\int_0^{T_j}g(R^j(t),\Phi^j(t))\,dt$,
where $T_j$ is the period of $\Gamma_j=(R^j,\Phi^j,\Psi^j)$. So, we require $b_1b_2<0$ (in fact, only $b_1\ne b_2$
is enough, as we show). Then, our orbit selection method works as follows: assume that
initially the $j$-th soliton is in the interval $[L_j^-,L_j^+]$ with $L_j^+-L_j^-$ large enough;
then until it remains in that interval, we allow the internal state $(R,\Phi,\Psi)$ of the soliton to
jump randomly between $\Gamma_1$ and $\Gamma_2$ along the heteroclinic orbits
(thus we gain the complexity which is enough to have the positive entropy); however, when the soliton reaches the bound
(say, $L_j^+$), we stop allowing jumps and consider only orbits that stay near $\Gamma_1$
until the soliton position $p_j(t)$ arrives close to $(L_i^++L_i^-)/2$ (when the bound $L_j^-$ is achieved, the orbit
must stay near $\Gamma_2$); after $p_j(t)$ is driven to the middle of the interval,
the random motion is allowed again, and so on.

We proved in Section \ref{s3} that the above described procedure can be implemented {\it simultaneously} for
all chaotic solitons on the grid, and it
allows indeed for a selection of a set of spatially non-walking solitons with positive space-time entropy.
In order to do this, we need a further development of the theorem on normally-hyperbolic manifolds in
the countable product of partially hyperbolic sets which is proved in Section \ref{s2}; namely we prove
certain ``asymptotic phase'' results in Section \ref{s2.5}.

As the above discussion shows, the theory we build is readily applicable to any dissipative PDE
for which the weak soliton interaction system for some finite multi-soliton configuration exhibits
a chaotic dynamics. In analogy to the finite-dimensional case, we are now able to
analyze localized structures and effectively use them for the understanding of
space-time dynamics generated by PDEs.

We are grateful to A.Mielke and A.Vladimirov for useful discussions,
and to WIAS (Berlin) and BenGurion University for the hospitality.

\section{Space-time chaos in complex Ginzburg-Landau equation with broken phase symmetry}\label{s4}
Consider the one-dimensional complex Ginzburg-Landau equation
\begin{equation}\label{4.ppcGL}
\Dt u=(1+i\beta)\partial_{x}^2 u-(1+i\delta)u-uH(|u|^2)+\mu G(u)
\end{equation}
where $u=u_1+iu_2$ is an unknown function of $x\in\R$ and $t\in\R$, the function $H:\R\to\Bbb C$
is smooth, $H(0)=0$,
and parameters $\beta$, $\delta$ are real; the symmetry-breaking parameter $\mu$
is assumed to be small, and the function $G$ is smooth.

Let $\Cal K$ be a set of solutions of \eqref{4.ppcGL} which are defined
and uniformly bounded for all $(t,x)\in\R^2$
(under certain standard dissipativity assumptions, equation \eqref{4.ppcGL} will have global attractor; in
this case one can choose as the set $\Cal K$ the set of all solutions that lie in the attractor,
see more after Theorem \ref{Th05.entropy}). The complexity of spatio-temporal behavior of the solutions
can be characterized by the space-time topological entropy defined as
\begin{equation}\label{4.ent}
h_{s-t}(\Cal K)=
\lim_{\eb\to0}\limsup_{(T,R)\to\infty}\frac1{4TR} h_\eb(\Cal K\big|_{|t|\le T,|x|\le R})
\end{equation}
where $\Cal K\big|_{|t|\le T,|x|\le R}$ stands for the set of functions from $\Cal K$ restricted on
the space-time window $\{|t|\leq T, |x|\leq R\}$, and $h_\eb$ denotes the Kolmogorov
$\varepsilon$-entropy of this set, i.e. the logarithm of
a minimal number of $\eb$-balls in the space $L^\infty([-T,T]\times [-R,R])$ which are necessary to cover
\footnote{it follows in a standard way from the parabolic regularity, that $\Cal K\big|_{|t|\le T,|x|\le R}$
is compact} the set; see \cite{255}. It is well-known (see, e.g., \cite{15,281,55,56})
that the space-time topological entropy $h_{s-t}(\Cal K)$ is well-defined and {\it finite} in our case.
Thus, if $h_{s-t}(\Cal K)$ is strictly positive for some set $\Cal K$,
then the number of various spatio-temporal patterns that are supported
by the equation grows {\em exponentially} with the volume of the space-time window.

In our construction of spatio-temporal chaos we assume that the nonlinearity $H$
is such that for some $\beta=\beta_0$ and $\delta=\delta_0$ the Ginzburg-Landau
equation \eqref{4.ppcGL} possesses at $\mu=0$ a stationary, spatially localized solution $u=U(x)$:
\begin{equation}\label{4.cGLv}
(1+i\beta_0)\partial_{x}^2 U-(1+i\delta_0)U-UH(|U|^2)=0;
\end{equation}
for the existence results see \cite{AfM1,AfM2} and references therein, and Theorem \ref{thnls}.

Equation \eqref{4.ppcGL} is invariant with respect to spatial translations $x\to x-\xi$ and, at $\mu=0$,
with respect to phase shifts $u\to e^{i\phi}u$. So, along with the given soliton $U(x)$, equation
\eqref{4.ppcGL} possesses at $\mu=0$ a family of stationary solitons:
\begin{equation}\label{4.man1}
u=U_{\xi,\phi}(x):=e^{i\phi}U(x-\xi), \ (\xi,\phi)\in\R^1\times\Bbb S^1.
\end{equation}
Because of the symmetry with respect to $x\to -x$, along with
the soliton $u=U(x)$, equation \eqref{4.cGLv} also has a localized solution $u=U(-x)$.
Equation \eqref{4.cGLv} is an ODE with 4-dimensional phase space. A localized solution
corresponds to a homoclinic intersection of the stable and unstable manifolds of the zero
equilibrium of this system. Since these manifolds are 2-dimensional and family \eqref{4.man1} is
2-parametric, all the localized solutions of \eqref{4.cGLv} are contained in family \eqref{4.man1}.
Thus, $U(-x)\equiv e^{i\phi_0}U(x-\xi_0)$ for some $\phi_0,\xi_0$, which immediately implies
that $U_{\xi/2,0}(-x)=\pm U_{\xi/2,0}(x)$. In other words, we may from the very beginning assume
that our soliton is chosen such that it is either symmetric:
\begin{equation}\label{4.syms}
U(-x)\equiv U(x),
\end{equation}
or antisymmetric ($U(-x)\equiv -U(x)$). In this paper we consider the symmetric case, i.e. we assume that
\eqref{4.syms} holds (in the antisymmetric case the soliton interaction equations are different; however
one can show that a small perturbation of an equation with antisymmetric soliton creates symmetric solitons -
cf. \cite{AfM1,AfM2}, so the results of our paper can be applied in this way).

Since every function in \eqref{4.man1} is a stationary solution of \eqref{4.cGLv} at $\mu=0$,
it follows that the functions $\varphi_1:=-\partial_{\xi}U_{\xi,0}\big|_{\xi=0}=\partial_x U$ and
$\varphi_2:=\partial_{\phi}U_{0,\phi}\big|_{\phi=0}=i U$ belong to
the kernel of the linearization $\Cal L_U$ of \eqref{4.cGLv} at $U$:
$\displaystyle\;\;\Cal L_{_U}\varphi_{1,2}=0,\;\;$ where
\begin{equation}\label{lulin}
\Cal L_{_U}\varphi:=(1+i\beta_0)\partial_{x}^2\varphi-(1+i\delta_0)\varphi-H(|U|^2)\varphi-
|U|^2H'(|U|^2)\varphi-U^2H'(|U|^2)\bar \varphi
\end{equation}
($\bar \varphi$ is a complex conjugate to $\varphi$). Thus, zero is a double eigenvalue of $\Cal L_U$.

We assume that the soliton $U$ is {\em non-degenerate} in the sense that the rest of the spectrum of
$\Cal L_U$ is bounded away from the imaginary axis; e.g. the algebraic multiplicity
of the zero eigenvalue is two
(note that since $U(x)\to0$ as $x\to\pm\infty$, the operator $L_U$
is a compact perturbation of the operator
$\varphi\mapsto(1+i\beta_0)\partial_{x}^2\varphi-(1+i\delta_0)\varphi$,
so the essential spectrum is bounded away
from the imaginary axis; however, one should check that the eigenvalues stay away
from the imaginary axis as well).

Under the non-degeneracy assumption, the conjugate operator $\Cal L^{\dagger}_U$, which we define as
\begin{equation}\label{lulc}
\Cal L^{\dagger}_U\psi:=
(1+i\beta_0)\partial_{x}^2\psi-(1+i\delta_0)\psi-H(|U|^2)\psi-|U|^2H'(|U|^2)\psi-\bar U^2H'(|U|^2)\bar \psi,
\end{equation}
also has a two-dimensional kernel. The corresponding pair of adjoint
eigenfunctions $\psi_1$ and $\psi_2$ can be chosen such that
\begin{equation}\label{normaf}
(\varphi_i,\psi_j):=\Ree \int_{-\infty}^{+\infty} \varphi_i(x)\psi_j(x)dx=\delta_{ij},\ \
\psi_1(-x)=-\psi_1(x),\ \ \psi_2(-x)=\psi_2(x).
\end{equation}

As $x\to\pm\infty$, the functions $U$, $\varphi_i$, $\psi_i$ decay exponentially,
with the rate $\lambda$ given by
\begin{equation}\label{4.quadra}
\Ree\lambda=-\alpha<0,\;\;\Imm\lambda=\omega,\;\;\;\; (-\alpha+i\omega)^2(1+i\beta_0)=(1+i\delta_0),
\end{equation}
see \cite{AfM1,57} for details. Thus, we have
\begin{equation}\label{4.tails}
U\sim re^{(-\alpha+i\omega)|x|},\;\; \psi_1\sim se^{(-\alpha+i\omega)|x|}\mbox{ sign}(x),\;\;
\psi_2\sim qe^{(-\alpha+i\omega)|x|}\mbox{   as }|x|\to\infty,
\end{equation}
where $r,s,q$ are some non-zero complex constants. We introduce the notation
\begin{equation}\label{4.abt}
ae^{i\theta_1}:=4isr(1+i\beta_0)\lambda,\ \ be^{i\theta_2}:=4iqr(1+i\beta_0)\lambda,\ \ \theta:=\theta_2-\theta_1.
\end{equation}

Denote
\begin{equation}\label{4.ccf}
F(\phi):=\Ree \int_{-\infty}^{+\infty}e^{-i\phi}\psi_{2}(x) G(e^{i\phi}U(x))dx,
\end{equation}
where $G(u)$ is the symmetry-breaking term in \eqref{4.ppcGL}. Since $F(\phi)$ is periodic, the equation
\begin{equation}\label{psist}
F'(\phi^*+\frac{\pi}{4})+F'(\phi^*-\frac{\pi}{4})=0
\end{equation}
always has solutions. We assume that there is a solution $\phi^*$ such that
\begin{equation}\label{ccf}
c:=2F'(\phi^*+\frac{\pi}{4})\neq 0,
\end{equation}
\begin{equation}\label{ppf}
F''(\phi^*+\frac{\pi}{4})+F''(\phi^*-\frac{\pi}{4})\neq 0.
\end{equation}
Conditions \eqref{psist}-\eqref{ppf} define the constant $\phi^*$. Denote also
\begin{equation}\label{ccgm0}
\gamma:=\frac{1}{c}[F(\phi^*+\frac{\pi}{4})-F(\phi^*-\frac{\pi}{4})].
\end{equation}
In the basic case $G(u)\equiv 1$, we have $F(\phi)=\tilde c\cos(\phi-\zeta)$, where
$\tilde ce^{i\zeta}=\int_{-\infty}^{+\infty}\psi_2(x)dx$. It is easy to see that $\phi^*=\zeta$,
$c=-\tilde c \sqrt{2}$ and $\gamma=0$ in this case, and that both conditions \eqref{ccf} and \eqref{ppf}
are fulfilled provided $\left|\int_{-\infty}^{+\infty}\psi_2(x)dx\right|\neq 0$.

\begin{theorem}\label{Th05.entropy} Let, along with \eqref{ccf},\eqref{ppf}, the following conditions
be satisfied for a non-degenerate, symmetric stationary soliton $U(x)$:
\begin{equation}\label{rhocon}
a\neq0,\; b\neq0,\;\omega\neq 0,\; \omega\neq 2\gamma\alpha,\; \cos\theta\neq 0,\; \alpha\sin\theta+\omega\cos\theta\neq0,
\end{equation}
\begin{equation}\label{rhocon4}
4\omega\frac{a}{b}(\cos\theta+2\gamma\sin\theta)<\left[1+2\gamma\frac{a}{b}(\alpha\cos\theta+\omega\sin\theta\right]^2.
\end{equation}
Then, arbitrarily close to $\mu=0$ and $\delta=\delta_0$ there exist an interval of values of
$\mu$ and an interval of values of $\delta$ such that the corresponding equation \eqref{4.ppcGL}
has a uniformly bounded set of globally defined solutions with strictly positive space-time entropy.
\end{theorem}

\begin{proof} Each of the solutions of equation \eqref{4.ppcGL} that belong to the large (of positive
entropy) set we are going to construct can be viewed as a slowly evolving multi-soliton configuration.
Namely, we choose a sufficiently large $L$ and consider solutions $u(x,t)$ which for every $t\in\R$ stay close,
in the space $C_b(\R)$ of bounded continuous functions of $x$, to the {\em multi-soliton manifold} $\Bbb M_L$
defined as the set of all functions $u(x)$ of the form
\begin{equation}\label{4.mp}
u(x)=u_{\mm}:=\sum U_{\xi_j,\phi_j}:=\sum_{j\in\Bbb Z}e^{i\phi_j}U(x-\xi_j),
\end{equation}
where $\mm:=\{\xi_j,\phi_j\}_{j=-\infty}^{j=+\infty}$ is any sequence such that
\begin{equation}\label{4.sep}
\inf_{j\in\Bbb Z}(\xi_{j+1}-\xi_j)>2L.
\end{equation}
For sufficiently large $L$, the multi-soliton manifold is indeed an
infinite-dimensional submanifold of $C_b(\R)$ which is parameterized by the sequences
$\mm:=\{\xi_j,\phi_j\}$ of the soliton positions and phases (see \cite{57}). The boundary $\partial\Bbb M_L$
is given by $\inf_{j\in\Bbb Z}(\xi_{j+1}-\xi_j)=2L$.

We will seek for solutions of equation \eqref{4.ppcGL} in the form
$u(t):=u_{\mm(t)}+w(t)$
where $\mm(t)$ is a slow trajectory in $\Bbb M_L$ and $w(t)$ is a small corrector.
Recall a result from \cite{57}.

\begin{theorem}\label{Th4.main} For all $L$ large enough
there exists a $C^k$-map $\Bbb S:\Bbb M_L\to C_b(\R)$ such that
\begin{equation}\label{4.small}
\|\Bbb S\|_{C^k(\Bbb M_L,C_b(\R))}\le C e^{-\alpha L}
\end{equation}
(where $\alpha>0$ is defined by \eqref{4.quadra}) and that the manifold
$\Cal S:=\{u=u_{\mm}+\Bbb S(u_{\mm}),\; \mm\in\Bbb M_L\}$
is invariant with respect to equation \eqref{4.ppcGL}. Namely, there exists
a $C^k$-vector field $\Cal F$ on $\Bbb M_L$ such that given any solution of
\begin{equation}\label{4.rde}
\frac{d}{dt}\mm(t)={\Cal F}(\mm(t))
\end{equation}
defined on a time interval $t\in(t_-,t_+)$, the function
\begin{equation}\label{4.lift}
u_{\mm(t)}+\Bbb S(u_{\mm(t)}), \ \ t\in(t_-,t_+),
\end{equation}
solves equation \eqref{4.ppcGL}.

Moreover, system \eqref{4.rde} has the following form:
\begin{equation}
\label{4.lrde}
\begin{aligned}
\frac d{dt}\xi_j=2\Ree[\;s\;r\;(1+i\beta_0)\lambda \left\{e^{\lambda(\xi_{j+1}-\xi_j)+i(\phi_{j+1}-\phi_j)}-
e^{\lambda(\xi_j-\xi_{j-1})+i(\phi_{j-1}-\phi_j)}\right\}]+\dots,\\
\frac d{dt}\phi_j=-2\Ree[\;q\;r\;(1+i\beta_0)\lambda \left\{e^{\lambda(\xi_{j+1}-\xi_j)+i(\phi_{j+1}-\phi_j)}+
e^{\lambda(\xi_j-\xi_{j-1})+i(\phi_{j-1}-\phi_j)}\right\}]\qquad\\+\mu F(\phi_j)-(\delta-\delta_0)+\dots,
\end{aligned}
\end{equation}
where $\alpha$ and $\omega$ are the same as in \eqref{4.quadra}, the constants
$r,s,q$ are defined by \eqref{4.tails}, the function $F$ is defined by \eqref{4.ccf}, and the dots stand
for terms which are $O(e^{-3\alpha L}+\mu^2+(\delta-\delta_0)^2)$ in $C^k(\Bbb M_L,\R)$-metric,
uniformly for all $j\in\Bbb Z$.
\end{theorem}
The complete proof of formulas \eqref{4.lrde} occupies a substantial part of \cite{57}. For reader convenience,
we provide a brief heuristic derivation of the equations in Appendix.

According to Theorem \ref{Th4.main},
the evolution of well-separated multi-soliton configurations in the driven Ginzburg-Landau
equation is governed by system \eqref{4.lrde}. Therefore, in order to prove the positivity of space-time entropy
in equation \eqref{4.ppcGL}, it is enough to find a large set of solutions of system \eqref{4.lrde} which satisfy
the separation condition \eqref{4.sep}.\footnote{recall that system \eqref{4.lrde} is defined on the
manifold $\Bbb M_L$ whose boundary is given by \eqref{4.sep}; outside this boundary the reduction
to the invariant manifold $\Cal S$ may fail -- the so-called strong soliton interaction,
soliton collisions, etc., may take place}

The corresponding theory for a class of lattice dynamical systems which includes system \eqref{4.lrde} is
built in Sections \ref{s2}-\ref{s3}. In particular, Theorem \ref{Th3.1} gives a general result on the existence
of a set $\tilde{\Cal K}$ of {\em non-walking} trajectories of a lattice dynamical system such that
$h_{s-t}(\tilde{\Cal K})>0$. In what follows we will show that a certain subsystem of \eqref{4.lrde}
indeed satisfies conditions of Theorem \ref{Th3.1}.

For any $v$ and any sufficiently large $L$ we may define a sequence $L_n$, $n\in\Bbb Z$, as follows:
\begin{equation}\label{5.L}
L_0=ve^{-\alpha L}t, \ \ L_{2n+1}=L_{2n}+4L, \ \ L_{2n+2}=L_{2n+1}+2L.
\end{equation}
We will look for pulse configurations which satisfy
$\xi_j(t)=L_j+\eta_j(t)$ where
\begin{equation}\label{5.separated}
|\eta_j(t)|\le C, \ \ t\in\R, \ j\in\Bbb Z,
\end{equation}
where the constant $C$ is independent of $L$, $j$ and $t$. In other words, we have a grid of weakly interacting
pulse pairs with the distance between the pulses in the pair of order $2L$ and the distance between
pairs of order $4L$. Assumption \eqref{5.separated} then means that we should ensure that this structure
is preserved for all $t$ although a uniform spatial drift of the whole grid is allowed
($ve^{-\alpha L}$ is the velocity of the drift).

Further, we introduce the scaling
$\displaystyle \tau:=t e^{-2\alpha L}, \ \ \ \Omega:=(\delta-\delta_0) e^{2\alpha L}, \ \ \nu:=\mu e^{2\alpha L}$.
We will consider a region of bounded $\Omega$ and $\nu$, which corresponds to $\mu$ and $\delta-\delta_0$
of order $O(e^{-\alpha L})$. We also assume $L=\frac{\pi n}{\omega}$, $n\in\Bbb N$.
Equations \eqref{4.lrde} recast as follows (see \eqref{4.abt}):
\begin{equation}\label{4.spe}
\begin{aligned}
\frac d{d\tau}\eta_{2j+1}=v-\;\frac{a}{2} e^{-\alpha R_j}
\sin(\omega R_j-\Phi_j+\theta_1)\;+\;O(e^{-\alpha L}),\qquad\qquad\qquad\\
\frac d{d\tau}\eta_{2j+2}=v+\;\frac{a}{2} e^{-\alpha R_j}
\sin(\omega R_j+\Phi_j+\theta_1)\;+\;O(e^{-\alpha L}),\;\;\qquad\qquad\qquad\\
\frac d{d\tau}\phi_{2j+1}=\frac{b}{2} e^{-\alpha R_j}
\sin(\omega R_j-\Phi_j+\theta_2)+\nu F(\phi_{2j+1})-\Omega\;+\;O(e^{-\alpha L}),\\
\frac d{d\tau}\phi_{2j+2}=\frac{b}{2} e^{-\alpha R_j}
\sin(\omega R_j+\Phi_j+\theta_2)+\nu F(\phi_{2j+2})-\Omega\;+\;O(e^{-\alpha L}),
\end{aligned}
\end{equation}
where we denote $R_j:=\eta_{2j+2}-\eta_{2j+1}$, $\Phi_j:=-(\phi_{2j+2}-\phi_{2j+1})$.
As we see, only interaction {\em inside} the soliton pairs gives a contribution
into the leading terms of equations \eqref{4.spe}: since the distance between pairs is, in our configuration, of
order $4L$, the leading term for the interaction between solitons from different pairs will be of order
$O(e^{-4\alpha L})$ in the non-rescaled time $t$, so after the time rescaling it is of order $O(e^{-2\alpha L})$,
i.e. it is absorbed in the $O(e^{-\alpha L})$-terms in \eqref{4.spe}.

Let us rewrite the system in the coordinates $R_j$, $\Phi_j$, $\Psi_j:=(\phi_{2j+1}+\phi_{2j+2})/2$, and
$p_j:=(\eta_{2j+1}+\eta_{2j+2})/2$ (i.e. $p_j$ is the center of the soliton pair, $R_j$ is the distance between the solitons
in the pair, $\Phi_j$ and $\Psi_j$ describe the soliton phases). We obtain
\begin{equation}\label{4.finp}
\frac d{d\tau}p_j=v+\frac{a}{2}e^{-\alpha R_j}\cos(\omega R_j+\theta_1)\sin(\Phi_j)\;+\;O(e^{-\alpha L}),
\end{equation}
\begin{equation}\label{4.finy}
\left\{\!\!\!\!\begin{array}{l}\displaystyle
\frac {d R_j}{d\tau}=ae^{-\alpha R_j}\sin(\omega R_j\!+\!\theta_1)\cos(\Phi_j)\;+\;O(e^{-\alpha L}),\\
\\ \displaystyle
\frac{d\Phi_j}{d\tau}\!=\!be^{-\alpha R_j}\!\cos(\omega R_j\!+\!\theta_2)\!\sin(\Phi_j)\!
+\!\nu\!\left[F(\Psi_j\!+\!\frac{\Phi_j}{2})\!-\!F(\Psi_j\!-\!\frac{\Phi_j}{2})\right]\!+O(e^{-\alpha L}\!),\\
\\ \displaystyle
\frac{d\Psi_j}{d\tau}\!=\!\frac{b}{2} e^{-\alpha R_j}\!\sin(\omega R_j\!+\!\theta_2)\!\cos(\Phi_j)\!\!
+\!\frac{\nu}{2}\!\left[\!F(\Psi_j\!+\!\frac{\Phi_j}{2})\!+\!F(\Psi_j\!-
\!\frac{\Phi_j}{2})\!\right]\!\!-\!\Omega\!+O(e^{-\alpha L}\!),
\end{array}\right.
\end{equation}
At large $L$ this system is a lattice dynamical system of form \eqref{2.10}: at $L=+\infty$ the subsystems
that correspond to different $j$ are independent and identical, and the equations for variables
$\;y_j:=(R_j,\Phi_j,\Psi_j)\;$ (equations \eqref{4.finy}) are independent of the $p$-equation \eqref{4.finp}.
Therefore, in order to prove Theorem \ref{Th05.entropy}, it is enough to check that system \eqref{4.finp},\eqref{4.finy}
satisfies conditions of Theorem \ref{Th3.1} at some $v$. According to that theorem, we will then obtain,
for all sufficiently large $L$, the existence of a set $\tilde{\Cal K}$
of solutions of system
\eqref{4.finp},\eqref{4.finy} which has a positive space-time entropy and is uniformly bounded by a constant
independent of $L$ (i.e. condition \eqref{5.separated} is fulfilled -- this, in turn, ensures that the separation
condition \eqref{4.sep} holds, with somewhat smaller $L$, for all the solutions from $\tilde{\Cal K}$).
Now, lifting the set $\tilde{\Cal K}$ by formula \eqref{4.lift}, we obtain a uniformly bounded set $\Cal K$ of
globally defined solutions
of the perturbed Ginzburg-Landau equation, and the positivity
of the space-time entropy of the set $\Cal K$ follows from the smallness of $\Bbb S$
and the positivity of the space-time entropy of $\tilde{\Cal K}$.

Thus, to finish the proof we need the following
\begin{lemma}\label{Lemzero3}
Assume \eqref{psist}-\eqref{rhocon4}.
Then there exists an open region of values of $\nu$ and $\Omega$ for which the system
\begin{equation}\label{4.finj}
\frac{d}{d\tau}y:=\left\{\begin{array}{l}\displaystyle
\frac d{d\tau}R=ae^{-\alpha R}\sin(\omega R+\theta_1)\cos(\Phi),\\
\\ \displaystyle
\frac d{d\tau}\Phi=be^{-\alpha R}\cos(\omega R+\theta_2)\sin(\Phi)
+\nu\!\left[F(\Psi\!+\!\frac{\Phi}{2})\!-\!F(\Psi\!-\!\frac{\Phi}{2})\right],\\
\\ \displaystyle
\frac d{d\tau}\Psi=\frac{b}{2} e^{-\alpha R}\sin(\omega R+\theta_2)\cos(\Phi)
+\!\frac{\nu}{2}\!\left[\!F(\Psi\!+\!\frac{\Phi}{2})\!+\!F(\Psi\!-\!\frac{\Phi}{2})\!\right]\!\!-\!\Omega
\end{array}\right.
\end{equation}
behaves chaotically, i.e. it has a basic
(=non-trivial, uniformly-hyperbolic, compact, locally-maximal, transitive, invariant)
set $\Lambda$. Moreover, in $\Lambda$ one can find two periodic orbits, $y=y_-(\tau)$ and $y=y_+(\tau)$,
of periods $T_-$ and $T_+$, respectively, such that
\begin{equation}\label{5.g}
\frac{1}{T_-}\int_0^{T_-} g(y_-(\tau))d\tau\neq \frac{1}{T_+}\int_0^{T_+} g(y_+(\tau))d\tau,
\end{equation}
where $g(y):=v+\frac{a}{2}e^{-\alpha R}\cos(\omega R+\theta_1)\sin(\Phi)$.
\end{lemma}

One may check that condition \eqref{5.g} implies that
\begin{equation}\label{5.gv}
\int_0^{T_-} g(y_-(\tau))d\tau \; \cdot \; \int_0^{T_+} g(y_+(\tau))d\tau <0
\end{equation}
for an appropriately chosen $v$.
Hence, the lemma indeed establishes the required fulfilment of conditions of Theorem \ref{Th3.1}:
chaotic system \eqref{4.finj} coincides with the $y$-subsystem \eqref{4.finy} at $L=+\infty$ (for every $j$),
and condition \eqref{5.gv} coincides with condition \eqref{3.1dll} (the function $g$ is the right-hand side
of the $p$-equation \eqref{4.finp}). So, it remains to prove the lemma.

We note that numerically the existence of chaos in system \eqref{4.finj} with $F=c\cos\phi$ as well as different
scenarios of its emergence for various parameter values were established in \cite{46}.
In our {\em analytic} proof of chaotic
behavior we use one of the scenarios mentioned in \cite{46}.
Namely, we find an equilibrium of system \eqref{4.finj}
with 3 zero characteristic eigenvalues. It is known \cite{CTS,Ibanez}
that bifurcations of such equilibrium lead to a Shilnikov saddle-focus homoclinic loop, hence to chaos.

{\em Proof of Lemma \ref{Lemzero3}.}
For $\Omega=\frac{\nu}{2}[F(\phi^*+\frac{\pi}4)+F(\phi^*-\frac{\pi}4)]$ and $\nu$ such that
\begin{equation}\label{rn}
\cos Z=-\gamma\frac{c\nu}{b}e^{\alpha (Z-\theta_2)/\omega},
\end{equation}
system \eqref{4.finj} has an equilibrium state at $\Phi=\frac{\pi}{2}$, $\Psi=\phi^*$, $R=(Z-\theta_2)/\omega$,
where $\phi^*$ is given by \eqref{psist} (see also \eqref{ccf},\eqref{ccgm0})
By \eqref{psist},\eqref{ccf}, the linearization matrix at such equilibrium is
$$\left(\begin{array}{ccc}
0 & -\rho_1 & 0\\
-\rho_2 & 0 & c\nu \\
0 & -\frac{b}{2} e^{-\alpha (Z-\theta_2)/\omega}\sin Z+\frac{1}{4}c\nu & 0
\end{array}\right),$$
where $\rho_1\!:=\!a e^{-\alpha (Z-\theta_2)/\omega}(\sin Z\cos\theta\!-\!\cos Z\sin\theta)$,
$\rho_2\!:=\!b\omega e^{-\alpha (Z-\theta_2)/\omega}(\alpha\cos Z\!+\!\omega\sin Z)$.
This matrix has three zero eigenvalues at $\nu=\nu^*$ provided
\begin{equation}\label{d2minor}
D(\nu^*):=\rho_1\rho_2+c\nu^* (\frac{1}{4}c\nu^*-\frac{b}{2} e^{-\alpha(Z-\theta_2)/\omega}\sin Z)=0.
\end{equation}
At $\gamma\neq0$ system \eqref{d2minor},\eqref{rn} for $Z=Z^*$ transforms into
$$\cos^2Z^* [1-4\gamma^2\alpha\frac{a}{b}\sin\theta]+
2\gamma \sin Z^*\cos Z^* [1+2\gamma\frac{a}{b}(\alpha\cos\theta-\omega\sin\theta)]
+4\gamma^2\frac{a}{b}\omega\cos\theta \sin^2 Z^*=0,$$
and it is easy to check that the solvability of this equation is given by condition \eqref{rhocon4}.
Moreover, solutions satisfy
\begin{equation}\label{dn02}
D'(\nu^*)\neq 0.
\end{equation}
If $\gamma=0$, condition \eqref{rn} gives $\cos Z^*=0$, and one may check that condition \eqref{rhocon4}
in this case guarantees the solvability of equation \eqref{d2minor} for $\nu^*$ and the fulfillment of \eqref{dn02}.
It follows from \eqref{rhocon} that
$\rho_{1,2}\neq 0$ at the solutions (hence $\nu^*\neq0$) and that
\begin{equation}\label{shfcon}
\alpha\cos(Z^*-\theta)+\omega\sin(Z^*-\theta)\neq 0.
\end{equation}

At $\nu=\nu^*$ (the triple zero bifurcation moment) the vectors
$$v_1=\left(\begin{array}{c} -c\nu \rho_1 \\ 0 \\
-\rho_1\rho_2\end{array}\right),\quad
v_2=\left(\begin{array}{c} 0 \\  c\nu \\ 0\end{array}\right),\quad
v_3=\left(\begin{array}{c} 0 \\  0 \\ 1 \end{array}\right),$$
form a Jordan base. At $\nu$ close to $\nu^*$, take $Z$ satisfying \eqref{rn} and close to $Z^*$, and denote
\begin{equation}\label{uvw}
\left(\begin{array}{c} R-(Z-\theta_1)/\omega \\ \Phi-\frac{\pi}{2} \\ \Psi-\phi^*\end{array}\right)=
y_1 v_1+y_2v_2+y_3v_3=\left(\begin{array}{c} -c\nu \rho_1 \;y_1 \\
c\nu \;y_2\\
y_3 - \rho_1\rho_2 \;y_1 \end{array}\right).
\end{equation}
System (\ref{4.finj}) takes the form
\begin{equation}\begin{array}{l} \displaystyle
\dot y_1 = y_2 + O(y^2),\\ \displaystyle \dot y_2 = y_3 +O(y^2),\\ \displaystyle
\dot y_3  = \varepsilon_1
+ \varepsilon_2 y_2 +\rho y_1^2 + O(|y_1|^3+|y_1|(|y_2|+|y_3|)+y_2^2+y_3^2),
\end{array} \label{i4}
\end{equation}
where $\rho=\frac{1}{4} (\rho_1\rho_2)^2\nu [F''(\phi^*+\frac{\pi}4)+F''(\phi^*-\frac{\pi}4)]\neq 0$, and
$\varepsilon_1= \frac{1}{2}\nu [F(\phi^*+\frac{\pi}4)+F(\phi^*-\frac{\pi}4)]-\Omega$,
$\varepsilon_2=D(\nu)$,
i.e. $(\varepsilon_1,\varepsilon_2)$ are small parameters which are related by a diffeomorphism to
the original parameters $\nu$ and $\Omega$ near
the triple zero bifurcation moment (see \eqref{d2minor},\eqref{dn02}).

Scale the parameters as follows:
\begin{equation}\label{tav}
\varepsilon_1=-\frac{1}{\rho}s^6,\quad \varepsilon_2=E s^2
\end{equation}
for a sufficiently small $s$, and for some bounded $E$. By scaling the time
and the variables:
$$\tau\rightarrow \sigma/s,\qquad y_1\rightarrow  Y\varepsilon_1/s^3, \quad y_2\rightarrow Y_2 \varepsilon_1/s^2,
\quad y_3\rightarrow Y_3\varepsilon_1/s,$$
we bring system (\ref{i4}) to the form
\begin{equation}\label{i5a}
Y^{\prime\prime\prime} = 1-Y^2 + E Y^{\prime} + O(s)
\end{equation}
(where ${}^\prime$ denotes the differentiation with respect to the new, slow time $\sigma$).

The limit equation
\begin{equation}\label{i5}
Y^{\prime\prime\prime} = 1-Y^2 + E Y^{\prime}
\end{equation}
has two hyperbolic equilibria: $O_+:\;Y=1$, with a one-dimensional
stable manifold $W^s_+$ and a two-dimensional unstable manifold $W^u_+$, and $O_-:\;Y=-1$, with a two-dimensional
stable manifold $W^s_-$ and a one-dimensional unstable manifold $W^u_-$. At $E<3$ these equilibria are
{\em saddle-foci}, i.e. each of them has a pair of complex characteristic exponents. By \cite{Tsuzuki},
equation \eqref{i5} has, at $E=E^*=-\frac{19}{\sqrt[3]{2475}}$, a solution
$$Y(t) = -\;\frac{9}{2} \tanh(\sqrt[3]{11/120}\;t) + \frac{11}{2} \tanh^3(\sqrt[3]{11/120}\;t)$$
which connects the saddle-focus $O_-$ with $O_+$. This solution
corresponds to a curve $\Gamma_{-+}$ along which the one-dimensional manifolds $W^u_-$ and $W^s_+$ coincide.
By \cite{Ibanez}, at the same $E$ there exists another heteroclinic curve, $\Gamma_{+-}$,
which corresponds to a {\em transverse} intersection of the two-dimensional manifolds $W^u_+$ and $W^s_-$.
By the transversality, the heteroclinicnic orbit $\Gamma_{+-}$ persists for all $E$ close to $E^*$.
The other heteroclinic orbit, $\Gamma_{-+}$, splits as $E$ varies, and this results \cite{B1,B2} in the
sequence of values $E_k\rightarrow E^*$ which correspond to the existence of homoclinic
loops to the saddle-foci $O_+$ and $O_-$ (equation (\ref{i5}) is time-reversible, so homoclinic loops
to the both saddle-foci appear simultaneously). One can view the one-parameter family (\ref{i5})
as a smooth curve in the space of smooth flows in $\R^3$; then the parameter values $E_k$ correspond to
the intersections of this curve with smooth codimension-one surfaces filled by systems with a homoclinic loop
to, say, the saddle-focus $O_-$. Importantly, these intersections are transverse.
Therefore, fixing any arbitrarily large $k$, we will have at some $E$ close to $E_k$
a homoclinic loop to a saddle-focus close to $O_-$, for any one-parameter family
which is sufficiently close to (\ref{i5}).

Thus, given any sufficiently large $k$, at $E=E_k+O(s)$ equation (\ref{i5a}) has,
for every sufficiently small $s$,
a homoclinic loop $\Gamma^k_s$ to the saddle-focus $O_-$ at $Y=Y_{-}(k,s)=-1+O(s)$.
Denote as $\xi_{1,2,3}$ the characteristic
exponents at the saddle-focus, $\xi_1>0$, $\Ree \xi_2 =\Ree \xi_3<0$, $\Imm \xi_2=-\Imm\xi_3\neq 0$.
As $\xi_1+\xi_2+\xi_3\approx0$ here (the limit equation (\ref{i5}) is volume-preserving),
the Shilnikov condition of chaos, $\xi_1+\Ree\xi_2>0$, is automatically fulfilled.
Hence, by \cite{Sh1,Sh2} we obtain
an open region in the parameter plane which corresponds to a chaotic behavior (i.e. to the sought basic
hyperbolic set $\Lambda$) in equation (\ref{i5a}) and, equivalently, in the original system \eqref{4.finj}.

To finish the proof we need to show that the set $\Lambda$ can be chosen in such a way
that it will contain a pair of periodic orbits for which \eqref{5.g} is satisfied.
According to Remark \ref{Rem3.1d}, it is enough to check that the integral
of the function $g-g|_{_{O^-}}$ along the homoclinic loop to the saddle-focus $O_-$ is non-zero. In order
to verify this condition, let us rewrite the function $g$ in the new variables $(Y,Y^{\prime},Y^{\prime\prime})$:
$$g(Y,Y^\prime,Y^{\prime\prime})=v+\frac{a}{2}e^{-\alpha (Z-\theta_2)/\omega}\cos(Z-\theta)+
C\;s^3\;Y+O(s^6),$$
where $C:=-\;\frac{ ac\nu\rho_1}{2\rho}e^{-\alpha (Z-\theta_2)/\omega}(\alpha\cos(Z-\theta)+\omega\sin(Z-\theta))\neq 0$
(see (\ref{shfcon})).
Let $Y=Y(\sigma;k,s)$ be the solution of (\ref{i5a}) that corresponds to the homoclinic loop $\Gamma^k_s$;
note that $Y(\sigma;k,s)\to Y_-(k,s)$ exponentially as $\sigma\to\pm\infty$. Note also that
$Y(\sigma;k,s)=Y(\sigma;k,0)+O(s)$, therefore $\displaystyle
\int_{-\infty}^{+\infty}\left[g(Y(\sigma;k,s),Y^\prime(\sigma;k,s),Y^{\prime\prime}(\sigma;k,s))-
g(Y_-(k,s),0,0)\right]d\sigma\;\;=$\\ $\displaystyle =Cs^3\!\int_{-\infty}^{+\infty}\!\!\!(Y(\sigma;k,0)+1)d\sigma+O(s^4)$.
As $k\rightarrow+\infty$, the homoclinic loops of equation (\ref{i5}) approach
the heteroclinic cycle $\Gamma_{+-}\cup\Gamma_{-+}\cup O_-\cup O_+$ at $E=E^*$,
so the homoclinic loop $\Gamma^k_0$ to $O_-:\{Y=-1\}$ spends at large $k$
a large time in a neighborhood of the other equilibrium, $O_+:\{Y=+1\}$. Therefore, the integral of
$(Y(\sigma;k,0)+1)$ tends to $+\infty$ as $k\to+\infty$. Thus, for sufficiently large $k$ and sufficiently small $s$,
\begin{equation}\label{gloop}
\int_{-\infty}^{+\infty}
\left[g(Y(\sigma;k,s),Y^\prime(\sigma;k,s),Y^{\prime\prime}(\sigma;k,s))-g(Y_-(k,s),0,0)\right]d\sigma\neq 0.
\end{equation}
By Remark \ref{Rem3.1d}, this proves the lemma,
which finishes the proof of the theorem as well.
\end{proof}

The proof of the following proposition is standard, cf. \cite{54}.
\begin{proposition}\label{Prop5.attr} Let the non-linearity $H$ satisfy
\begin{equation}\label{5.dg}
\Ree H(z)\cdot z\ge -C; \ \ \ \ \ \ |H(z)|\le C(1+z^2),\ \ z\in \R_+
\end{equation}
for some constant $C$ independent of $z$. Then for all sufficiently small $\mu$ equation
\eqref{4.ppcGL} is well-posed in the space $C_b(\R)$ of uniformly bounded continuous functions and
generates a dissipative semigroup $S(t)_{t\ge0}$ in $C_b(\R)$, and this
semigroup possesses a global attractor $\Cal A$.
\end{proposition}
The attractor is defined here as follows.
Let $S(t)$, $t\geq0$, be a semigroup acting on the space $C_b(\R)$. A set $\Cal A\subset C_b(\R)$
is a global (locally-compact) attractor of this semigroup if
\par
1) $\Cal A$ is bounded in $C_b(\R)$ and compact in $C_{loc}(\R)$;
\par
2) $\Cal A$ is strictly invariant: $S(t)\Cal A=\Cal A$, $t\ge0$;
\par
3) as $t\to\infty$, the set $\Cal A$ attracts, in the topology of $C_{loc}(\R)$,
the images of all bounded subsets $B\subset C_b(\R)$, i.e. for
every neighborhood $\Cal O$ of $\Cal A$ in the local
topology and for every bounded $B\subset C_b(\R)$ there is a time $T=T(\Cal O, B)$ such that
$S(t)B\subset\Cal O(\Cal A)$ for all $t\ge T$.
\begin{remark} It is well-known (see e.g. \cite{281,55}) that, in contrast to the case of bounded domains,
the global attractor is usually not compact in $C_b(\R)$ if the underlying domain is unbounded. However,
attractor's restrictions to every bounded subdomain remain compact. The
attraction property itself holds, too, in this local topology only.
\end{remark}

A characteristic property of the global attractor is that it {\em consists of all initial conditions
which give rise to globally defined solutions}. Namely, a function $u_0(x)\in C_b(\R)$ belongs
to the attractor if and only if there exists a function $u(t,x)\in\Cal K$ such that $u_0(x)\equiv u(0,x)$.
Note that due to the invariance of the equation with respect to temporal and spatial translations,
the boundedness and local compactness of the attractor mean also that the set $\Cal K$ of the solutions which
are defined and bounded for all $(t,x)\in\R^2$ is bounded in $C_b(\R^2)$ and compact in $C_{loc}(\R^2)$.

Thus, we may define the space-time entropy of the attractor as the space-time entropy of the set $\Cal K$:
$h_{s-t}(\Cal A):=h_{s-t}(\Cal K)$ (see \eqref{4.ent}; more discussion and a comparison with other
definitions can be found e.g. in \cite{281,55}). As we mentioned (see \cite{15,56}),
the space-time entropy of the attractor of the Ginzburg-Landau equation is finite:
\begin{equation}\label{5.tf}
h_{s-t}(\Cal A) < \infty.
\end{equation}
The next Section gives an explicit example of a scientifically relevant equation with
\begin{equation}\label{5.tff}
h_{s-t}(\Cal A) > 0.
\end{equation}
\section{Attractor of positive space-time entropy in a perturbed nonlinear Shr\"odinger equation}\label{nls}
Here we prove the following
\begin{theorem}\label{thnls} Given any sufficiently large $\beta$, there exist (continuously depending on $\beta$)
intervals of values of $\delta$, $\rho$, $\varepsilon_1>0$, $\varepsilon_2$ and $\mu$ such that the attractor of
the equation
\begin{equation}\label{nlsgl}
\Dt u=(1+i\beta)\partial_{x}^2 u-(1+i\delta)u+(i+\rho)|u|^2u-(\varepsilon_1+i\varepsilon_2)|u|^4u+\mu
\end{equation}
has strictly positive space-time entropy.
\end{theorem}

\begin{proof} The global attractor of equation \eqref{nlsgl} exists at $\varepsilon_1>0$ according to Proposition
\ref{Prop5.attr}. By theorem \ref{Th05.entropy}, in order to prove \eqref{5.tff}
it is enough to show that the equation
\begin{equation}\label{nlsglstat}
(1+i\beta)\partial_{x}^2 U-(1+i\delta)U+(i+\rho)|U|^2U-(\varepsilon_1+i\varepsilon_2)|U|^4U=0
\end{equation}
has a non-degenerate symmetric localized solution at some $\delta$ that depends on the other parameters
$\beta$, $\rho_{1,2}$, $\varepsilon_{1,2}$, and that conditions \eqref{ccf},\eqref{rhocon},\eqref{rhocon4} are satisfied
at $\gamma=0$. The localized
solution of the ODE \eqref{nlsglstat} corresponds to an intersection of the two-dimensional stable and unstable
manifolds of the hyperbolic equilibrium at $U=0$. Because of the phase-shift symmetry, when these manifolds
intersect they coincide. The soliton non-degeneracy conditions imply (among other things) that as $\delta$ changes
the manifolds split {\em with a non-zero velocity}. It follows that a non-degenerate soliton will persist
at small perturbation of the nonlinearity, provided a small adjustment to the value of $\delta$ is made
(see more in \cite{AfM1,AfM2}). Thus, it is enough to consider the cubic equation
\begin{equation}\label{nlsc}
(1+i\beta)\partial_{x}^2 U-(1+i\delta)U+(i+\rho)|U|^2U=0;
\end{equation}
once the existence of a non-degenerate soliton is established for this equation, it can be carried on
to the equation \eqref{nlsglstat} for all sufficiently small $\varepsilon_{1,2}$,
and since conditions \eqref{ccf},\eqref{rhocon},\eqref{rhocon4} are open, they will persist as well.

Let us choose
$\displaystyle \beta=\frac{1}{B},\;\; \rho=\frac{B(1-2w^2)-3w}{1-2w^2+3wB},\;\;\delta=\frac{1-w^2+2wB}{B(1-w^2)-2w}$
for some small $B>0$ and $w$ such that $B>\frac{2w}{1-w^2}$.
Then, if we define
\begin{equation}\label{glnlssc}
\Cal U(x):=d_1U(x d_2)
\end{equation}
where $d_1=\sqrt{\frac{B(1-w^2)-2w}{1-2w^2+3wB}}$, $d_2=\sqrt{1-w^2-2\frac{w}{B}}$,
we obtain the following equation:
\begin{equation}\label{GL}
(i+B)\left[\partial_{x}^2 \Cal U - (1+i\omega)^2\Cal U + (1+i\omega)(2+i\omega) |\Cal U|^2\Cal U\right]=0.
\end{equation}
It has a localized stationary solution (see e.g. \cite{AfM1,AfM2})
\begin{equation}\label{sol}
U_*=\frac{1}{({\rm ch} (x))^{1+i\omega}}.
\end{equation}
The linearization operator $\Cal L$ is given by
\begin{equation}\label{lingl}
\begin{aligned}
\Cal L\varphi:=
(i+B)\left[\partial_{x}^2\varphi - (1+i\omega)^2\varphi + 2(1+i\omega)(2+i\omega) |U_*|^2\varphi+
\right.\qquad\qquad\qquad\\
\left.+\;(1+i\omega)(2+i\omega)^2\; U_*^2\;\bar\varphi\right].\quad
\end{aligned}
\end{equation}
The localized functions
\begin{equation}\label{evenodd}
\varphi_1(x)=-\partial_{x} U_*(x), \qquad \varphi_2(x)=iU_*(x)
\end{equation}
(the odd and, respectively, the even one) belong to the kernel of $\Cal L$.

We introduce a scalar product as
$\displaystyle \;\;(\varphi,\psi)={\Ree}\int_{-\infty}^{+\infty} \varphi(x)\psi(x) dx,\;\;$
so the conjugate to (\ref{lingl}) operator is
\begin{equation}\label{clingl}
\begin{aligned}
\Cal L^{\dagger}\psi:=(i+B)\left[\partial_{x}^2\psi - (1+i\omega)^2\psi +
2(1+i\omega)(2+i\omega) |U_*|^2\psi\right]+\qquad\qquad\qquad\\
+(-i+B)(1-i\omega)(2-i\omega) (\bar U_*)^2 \bar \psi.\quad
\end{aligned}
\end{equation}
As $\Cal L$ has two zero modes, one even and one odd,
the same holds true for the conjugate operator $\Cal L^{\dagger}$.
At $w=B=0$ the equation for zero eigenfunctions of $\Cal L^{\dagger}$ reads as
\begin{equation}\label{clinnls}
\partial_{x}^2\psi - \psi + 4 \Gamma^2\psi - 2 \Gamma^2 \bar \psi=0,
\end{equation}
where we denote
\begin{equation}\label{cosh}
\Gamma(x)=\frac{1}{{\rm ch} (x)};
\end{equation}
note that
\begin{equation}\label{23}
\Gamma''(x)=\Gamma-2\Gamma^3,\qquad \Gamma'''(x)=(1-6\Gamma^2)\Gamma'(x).
\end{equation}
It is easy to see that the odd and even localized solutions of \eqref{clinnls} are
\begin{equation}\label{evenoddnls}
\psi_1(x)= i\Gamma'(x), \qquad \psi_2(x)=\Gamma(x).
\end{equation}

We will look for asymptotic expansions of these solutions at small $w$
and $B$. By (\ref{clingl}), the localized zero modes of $\Cal L^\dagger$ satisfy
\begin{equation}\label{clingexp}
\begin{aligned}
\psi''(x) - \psi + 4\Gamma^2\psi - 2\Gamma^2 \bar\psi=
\qquad\qquad\qquad\qquad\qquad\qquad\qquad\qquad\qquad\qquad\\=i\omega\left[(2-6\Gamma^2)\psi
-(3\Gamma^2 + 4\Gamma^2\ln \Gamma) \bar\psi\right]+4iB \Gamma^2\bar\psi + O(w^2+B^2)
\end{aligned}
\end{equation}
(we take into account that $U_*$ depends on $w$ as well: by (\ref{sol}),(\ref{cosh})
$\bar U_*^2=\Gamma^2(1-2iw\ln \Gamma +O(w^2))$, while $|U_*|^2=\Gamma^2$). By (\ref{clingexp}), we have
\begin{equation}\label{psiuv}
\psi=u+iv+O(w^2+B^2),
\end{equation}
where
\begin{equation}\label{uveq}
\left\{\begin{array}{l}\displaystyle
u''(x) - u + 2\Gamma^2 u =-w v (2-3\Gamma^2 + 4\Gamma^2\ln \Gamma) + 4B v \Gamma^2, \\
v''(x) - v + 6\Gamma^2 v = w u (2-9\Gamma^2 - 4\Gamma^2\ln \Gamma) + 4B u \Gamma^2.
\end{array}\right.
\end{equation}
By (\ref{psiuv}),(\ref{uveq}),(\ref{23}) the two sought localized solutions of
(\ref{clingexp}) are given by
\begin{equation}\label{psievodd}
\psi_1=i\Gamma'(x)+S(x)+O(w^2+B^2),\qquad \psi_2=\Gamma(x)+iQ(x)+O(w^2+B^2),
\end{equation}
where $S$ and $Q$ are real, decaying to zero, as $x\rightarrow\pm\infty$,
functions which satisfy
\begin{equation}\label{p}
S'' - S + 2\Gamma^2 S =-w (2-3\Gamma^2 + 4\Gamma^2\ln \Gamma)\Gamma'(x) + 4B \Gamma^2 \Gamma'(x),
\end{equation}
\begin{equation}\label{q}
Q'' - Q + 6\Gamma^2 Q =w (2\Gamma-9\Gamma^3 - 4\Gamma^3\ln \Gamma) + 4B \Gamma^3.
\end{equation}
To find $S(x)$, we multiply (\ref{p}) to $\Gamma(x)$. The equation will take the form
(see (\ref{23})):
$$\Gamma S''-\Gamma''S=-w (\Gamma^2-\Gamma^4+\Gamma^4 \ln \Gamma)' +B (\Gamma^4)'.$$
By integrating this equation with respect to $x$, we find
$$\Gamma S'-\Gamma'S=-w (\Gamma^2-\Gamma^4+\Gamma^4 \ln \Gamma) +B \Gamma^4$$
(there is no integration constant in the right-hand side, since both $S$ and $\Gamma$
tend to zero as $x\rightarrow\pm\infty$). By solving the first-order equation, we
finally obtain
\begin{equation}\label{sfin}
\begin{array}{l}\displaystyle
S(x)=-w \Gamma(x) (x-\int \Gamma^2 dx +\int \Gamma^2 \ln \Gamma dx) +B \Gamma(x)\int \Gamma^2 dx\\
\displaystyle \qquad = -\frac{w}{{\rm ch}(x)} (2x-\frac{{\rm sh} (x)}{{\rm ch} (x)}(2-\ln{\rm ch} (x)))+
B\frac{{\rm sh} (x)}{{\rm ch}^2(x)}.
\end{array}
\end{equation}

Similarly, by multiplying (\ref{q}) to $\Gamma'(x)$ and integrating the obtained equation,
we find, with the use of (\ref{23}), that $\displaystyle\;\;\;
\Gamma'Q' - \Gamma''Q =w (\Gamma \Gamma'' - \Gamma^4\ln \Gamma) + B \Gamma^4.\;\;\;$
The solution is
\begin{equation}\label{qfin}
\begin{array}{l}\displaystyle
Q(x)=w (x\Gamma'(x)-\Gamma - \Gamma'(x)\int \frac{\Gamma^4 \ln \Gamma}{(\Gamma')^2} dx) +
B \Gamma'(x)\int \frac{\Gamma^4}{(\Gamma')^2} dx\\ \displaystyle
\qquad = -\;\frac{w}{{\rm ch}^2(x)} (2x {\rm sh} (x) + {\rm ch} (x) +
{\rm ch} (x) \ln{\rm ch} (x)) +B{{\rm ch} (x)}.\end{array}
\end{equation}
It is immediately seen that functions $S$ and $Q$ given by (\ref{sfin}),(\ref{qfin})
are localized indeed. Moreover, $S$ is odd and $Q$ is even, so by plugging
(\ref{sfin}) and (\ref{qfin}) in (\ref{psievodd}), we obtain the odd ($\psi_1$) and
even ($\psi_2$) zero eigenfunctions of $\Cal L^\dagger$.

One can also compute (see (\ref{evenodd})) that
\begin{equation}\label{sck1}
\begin{aligned}
\Ree \int_{-\infty}^{+\infty}\psi_1(x) \varphi_1(x) dx=\qquad\qquad\qquad\qquad\qquad\qquad\qquad\qquad\\
=-\int_{-\infty}^{+\infty} S\Gamma'dx+w\int_{-\infty}^{+\infty}(\Gamma')^2(1+\ln \Gamma) dx
+O(w^2+B^2)=\frac{2}{3}B+O(w^2+B^2),
\end{aligned}
\end{equation}
\begin{equation}\label{sck2}
\begin{aligned}
\Ree \int_{-\infty}^{+\infty} \psi_2(x) \varphi_2(x) dx=\qquad\qquad\qquad\qquad\qquad\qquad\qquad\qquad\\
=-w\int_{-\infty}^{+\infty} \Gamma^2\ln \Gamma dx-\int_{-\infty}^{+\infty}Q\Gamma dx +O(w^2+B^2)=
2(2w-B)+O(w^2+B^2).
\end{aligned}
\end{equation}
As we see, these inner products are non-zero for the values of $B$ and $w$ that we consider
here (small $B,w$ such that $B>0$ and $B>\frac{2w}{1-w^2}$). This
shows that there are no adjoint functions to the eigenfunctions (\ref{evenodd}). The absence (at small $B,w$)
of eigenvalues on the imaginary axis follows from \cite{225}. Thus, the pulse $\Cal U=U_*(x)$ is non-degenerate.

Returning to the non-rescaled variables,
we find that the soliton $U=d_1^{-1} U_*(x/d_2)$ of equation
\eqref{nlsc} is non-degenerate. The corresponding eigenfunctions of $\Cal L^\dagger_U$ are given by
$$\psi_1(x)=\frac{3d_1}{2B+O(w^2+B^2)}(i\Gamma'(x/d_2)+S(x/d_2)+O(w^2+B^2)),$$
$$\psi_2(x)=\frac{d_1}{2d_2((2w-B)+O(w^2+B^2))}(\Gamma(x/d_2)+iQ(x/d_2)+O(w^2+B^2))$$
(we normalize them so that (\ref{normaf}) is fulfilled, see (\ref{sck1}),(\ref{sck2})). By (\ref{sfin}),
(\ref{qfin}), we find
$$\psi_1(x)\sim-\;\frac{3d_1(i-w(2-\ln 2)-B+O(w^2+B^2))}{B+O(w^2+B^2)}e^{-(1+iw)|x|/d_2}\;{\rm sign}(x),$$
$$\psi_2(x)\sim\frac{d_1(1+iB-iw(1+\ln 2)+O(w^2+B^2))}{d_2(2w-B+O(w^2+B^2))}e^{-(1+iw)|x|/d_2}$$
as $x\rightarrow\pm\infty$, so $\omega=-w/d_2$, $\alpha=1/d_2$, and the coefficients $s$ and $q$
in \eqref{4.tails},\eqref{4.abt} are
\begin{multline*}
s=-\frac{3d_1(i-w(2-\ln 2)-B+O(w^2+B^2))}{B+O(w^2+B^2)}, \\
q=\frac{d_1d_2(1+iB-iw(1+\ln 2)+O(w^2+B^2))}{2w-B+O(w^2+B^2)}.
\end{multline*}
It is easy to see that all conditions \eqref{ccf},\eqref{rhocon},\eqref{rhocon4} hold at small
$w\neq 0$, $B>0$.
\end{proof}

\section{Normally-hyperbolic manifolds for lattice dynamical systems}\label{s2}

In this and the the next Sections we study a class of lattice dynamical systems which includes
systems describing weak interaction of solitons localized in space and chaotic in time,
e.g. system \eqref{4.finp},\eqref{4.finy}. We start with a skew-product system of ODE's
\begin{equation}\label{2.2}
y'(t)=f(y),\ \ p'(t)=g(y),
\end{equation}
where $f,g$ are $C^r$, $r\geq 1$. We assume that $y\in {\Bbb R}^n$, $p\in {\Bbb R}^m$;
for more clarity we will denote the space of $y$ variables as $Y$ and the space of $p$
variables as $P$. We will further assume that the $y$-part of our system:
\begin{equation}\label{2.1}
y'=f(y),
\end{equation}
possesses a bounded, uniformly-hyperbolic invariant set $\Lambda$.

Recall that the hyperbolicity means that for every point of $\Lambda$ there are
two subspaces, $N^s(y)$ and $N^u(y)$, such that the following holds:\\
1) $N^s(y)$ and $N^u(y)$ depend continuously on $y\in\Lambda$,\\
2) the direct sum of $N^s(y)$, $N^u(y)$ and $N^c(y):={\rm Span}(\dot y):=\{\lambda f(y)|\lambda\in\Bbb R\}$
constitutes the whole of $\Bbb R^n$,\\
3) given any orbit $y(t)$ from $\Lambda$, each of the families of subspaces $N^s(y(t))$ and $N^u(y(t))$
is invariant with respect to the flow of system \eqref{2.1} linearized about the orbit $y(t)$,\\
4) the linearized flow is exponentially contracting in restriction onto $N^s(y(t))$ as $t\rightarrow+\infty$
and in restriction onto $N^u(y(t))$ as $t\rightarrow-\infty$ (the flow, then, is expanding on
$N^s(y(t))$ as $t\rightarrow-\infty$ and on $N^u(y(t))$ as $t\rightarrow+\infty$).

The linearized system is
\begin{equation}\label{linvs}
\frac{d}{dt} v= f'(y(t))v.
\end{equation}
Since $v(t)=\dot y(t)=f(y(t))$ is a uniformly bounded solution of it, there exists a uniformly bounded
non-zero solution $y^*(t)$ for the conjugate system
\begin{equation}\label{linvc}
\frac{d}{dt} v = -f'(y(t))^\top v.
\end{equation}
As $y^*(t)$ solves \eqref{linvc}, it follows that $\frac{d}{dt} \<y^*(t)\cdot v(t)\> =0$
for every solution $v(t)$ of \eqref{linvs}, i.e. $\<y^*(t)\cdot v(t)\>$ stays constant.
Thus, since the solutions of \eqref{linvs} which lie
in $N^s(y(t))\oplus N^u(y(t))$ tend to zero either as $t\to+\infty$ or as $t\to-\infty$ and $y^*(t)$
is bounded, we find that this constant is zero for every $v\in N^s(y(t))\oplus N^u(y(t))$, i.e. the vector
$y^*(t)$ is orthogonal to $N^s(y(t))\oplus N^u(y(t))$ for all $t$. This condition defines $y^*$ up to
a scalar factor; we fix it by normalizing $y^*$ in such a way that
\begin{equation}\label{dotdot}
\<y^*(t)\cdot \dot y(t)\> \equiv 1.
\end{equation}

The exponential dichotomy for system \eqref{linvs} restricted to $v(t)\in N^s(y(t))\oplus N^u(y(t))$
implies that the equation
$$\frac{d}{dt} v(t)- f'(y(t))v= h(t)$$
has a unique uniformly bounded solution $v(t)\in N^s(y(t))\oplus N^u(y(t))$ for any uniformly
bounded function $h(t)\in N^s(y(t))\oplus N^u(y(t))$. It is more convenient for us to express this property
in the following equivalent way: the equation
\begin{equation}\label{basichyp}
\frac{d}{dt} v(t)- f'(y(t))v + \<y^*(t)\cdot v\> \dot y(t)= h(t)
\end{equation}
has a unique uniformly bounded solution $v(t)$ given any uniformly bounded function $h(t)$. More precisely,
equation \eqref{basichyp} defines a linear operator $L_{y}: h\mapsto v$ such that
\begin{equation}\label{clnrm}
\|v\|\leq C_\Lambda \|h\|.
\end{equation}
The assumed {\em uniform} hyperbolicity of the set $\Lambda$ means that the constant $C_\Lambda$ in
\eqref{clnrm} can be taken the same for all orbits $y\in\Lambda$.

\bigskip

Take a countable set of equations of type \eqref{2.2}. This produces
an uncoupled LDS ({\em lattice dynamical system}):
\begin{equation}\label{2.7}
y'_k(t)=f_k(y_k),\ \ p'_k(t)=g_k(y_k),\ \ k\in {\Bbb Z}
\end{equation}
We assume that the derivatives of $f_k$ and $g_k$ up to the order $r$ are uniformly continuous and bounded for
all $k$, and that for each $k$ the $k$-th individual ODE's in the LDS has a hyperbolic set $\Lambda_k$,
all these sets are uniformly bounded and uniformly
hyperbolic for all $k$ (the uniform hyperbolicity means in our approach that the constant $C_\Lambda$
in \eqref{clnrm} can be taken the same for all $k$).
In the example considered in Section \ref{s4}, the individual ODE's are identical to each other,
so the uniformity with respect to $k$ holds trivially.

By introducing Banach spaces
$$\begin{array}{l}
\Bbb Y:=l_\infty(Y),\ \ \|\yy\|_{\Bbb Y}:=\sup_{k\in\Bbb Z}\|y_k\|_Y,\ \ \yy:=\{y_k\}_{k\in\Bbb Z},\\
\Bbb P:=l_\infty(P),\ \ \|\pp\|_{\Bbb P}:=\sup_{k\in\Bbb Z}\|p_k\|_P,\ \ \pp:=\{p_k\}_{k\in\Bbb Z},
\end{array}$$
we may write the LDS as
$$\yy'(t)=\ff(\yy),\ \ \pp'(t)=\ggg(\yy),$$
where
$\ff:=\{f_k\}_{k\in\Bbb Z}, \ \ \ggg:=\{g_k\}_{k\in\Bbb Z}$.

The subject of our study will be a {\em coupled} LDS, obtained by a small smooth perturbation of
this system. Namely, we consider
\begin{equation}\label{2.10}
\begin{cases}
\yy'(t)=\ff(\yy)+\eb\FF_\eb(\yy,\pp),\\ \pp'(t)=\ggg(\yy)+\eb\GG_\eb(\yy,\pp),
\end{cases}
\end{equation}
where $\eb$ is a small parameter, and $\FF_\eb$ and $\GG_\eb$ are $C^r$-functions
$\Bbb Y\times\Bbb P\rightarrow \Bbb Y$ and, respectively, $\Bbb Y\times\Bbb P\rightarrow \Bbb P$;
by ``$C^r$'' we mean, here and below, that all the derivatives up to the order $r$ exist, are
uniformly continuous and uniformly bounded. We also assume continuity (in $C^r$) with respect
to $\varepsilon$, so
\begin{equation}\label{fgcr}
\|\FF_\eb\|_{C^r}+\|\GG_\eb\|_{C^r} \le C
\end{equation}
where $C$ is independent of $\eb$.

\medskip

Let $\yy^0(t):=\{y_k^0(t)\}_{k\in\Bbb Z}$ be a sequence of arbitrary orbits $y_k^0(t)\in\Lambda_k$; we will
say that $\yy^0(t)$ is an orbit from $\Lambda^\infty$. Each orbit $y_k^0(t)$ defines a curve in
the $Y$-space. The direct product of these curves, times the space $\Bbb P$, is a $C^r$-submanifold
of $\Bbb Y\times\Bbb P$, we will denote is as $\Bbb W^0_{\yy^0}$. Given an orbit $\yy^0$, the corresponding
manifold $\Bbb W^0_{\yy^0}$ is given by the equation
\begin{equation}\label{wman0}
y_k=y_k^0(\phi_k), \ k\in\Bbb Z,
\end{equation}
where the ``phases'' $\phi_k$ run all real values, independently for different $k$. If we introduce a
Banach space $\Psi$ of the bounded sequences $\Phi:=\{\phi_k\}_{k\in\Bbb Z}$ with the
uniform norm $\|\Phi\|:=\sup_{k\in\Bbb Z}|\phi_k|$,
then $\Bbb W^0_{\yy^0}$ is a $C^r$-embedding of $\Psi\times\Bbb P$ into $\Bbb Y\times\Bbb P$.
Obviously, $\Bbb W^0_{\yy^0}$ is invariant with respect to the non-coupled LDS \eqref{2.7}.
Moreover, this manifold is normally-hyperbolic (as each of the orbits $y^0_k$ is uniformly-hyperbolic).
It is a well-known general principle that normally-hyperbolic invariant manifolds persist
at small smooth perturbations (see \cite{Fenichel,HPS}).
The next theorem shows that this principle holds true in our setting.
\begin{theorem}\label{Th2.1} For all sufficiently small $\eb$, given
any orbit $\yy^0\in\Lambda^\infty$ there exists a uniquely defined $C^r$-manifold
${\Bbb W}_{\yy^0,\varepsilon}\subset {\Bbb Y}\times {\Bbb P}$, which is invariant
with respect to system \eqref{2.10}, depends continuously on $\eb$ (in $C^r$,
uniformly with respect to $\yy^0$), and coincides with ${\Bbb W}^0_{\yy^0}$ at $\varepsilon=0$.
Namely, ${\Bbb W}_{\yy^0,\varepsilon}$ is given by
\begin{equation}\label{wman}
y_k={\Bbb U}_k(\Phi,\pp,\varepsilon):=y_k^0(\phi_k)+{\Bbb V}_k(\Phi,\pp,\varepsilon),
\end{equation}
where
\begin{equation}\label{2.25}
\|{\Bbb V}_k\|_{C^{r-1}}=O(\varepsilon),\qquad \|{\Bbb V}_k\|_{C^{r}}=o(1)_{\eb\to 0},
\end{equation}
uniformly for all $k\in\Bbb Z$ and all $\yy^0\in\Lambda^\infty$.
\end{theorem}
\begin{proof}
We start with some preliminary constructions. Define {\em the exponential $\alpha$-norm}
$\|h\|_{\alpha}:=\sup_{t\in\Bbb R} e^{-\alpha|t|}\|h(t)\|$
on the space of continuous, uniformly bounded functions $h$; e.g. $\|\cdot\|_0$ is just the $C^0$-norm.
\begin{lemma}\label{l1m1} For all small $\alpha\geq 0$ and $\nu\geq 0$, for all functions
$\phi(t)$ such that
\begin{equation}\label{fnu}
|\phi'(t)-1|\leq\nu\quad \mbox{ for all  } t\in\Bbb R,
\end{equation}
and for any $A(t)$ and $b(t)$ sufficiently close (in $C^0$) to $f'(y(t))$ and, respectively,  to $y^*(t)$,
the equation
\begin{equation}\label{aphphi}
\frac{d}{dt} v(t)- A(\phi(t))v(t) + \<b(\phi(t))\cdot v(t)\> f(y(\phi(t))) = h(t)
\end{equation}
is uniquely solvable for any uniformly bounded function $h(t)$, and the corresponding linear operator
$L_{\phi}: h\mapsto v$ satisfies
\begin{equation}\label{clnre}
\|v\|_{\alpha}\leq C_\Lambda \|h\|_{\alpha}.
\end{equation}
Moreover, the operator $L_\phi$ is Lipshitz with respect to $\phi$:
if $v_1(t)$ and $v_2(t)$ are the solutions of equation \eqref{aphphi} which correspond to two different
functions $\phi_1(t)$ and $\phi_2(t)$ (and to the same right-hand side $h$), then
\begin{equation}\label{aplip}
\|v_2-v_1\|_{\alpha}\leq K \|h\|_0 \|\phi_2-\phi_1\|_{\alpha}
\end{equation}
for some constant $K$, proportional to the $C^1$-norms of $A$, $b$ and $f$.
\end{lemma}
\begin{proof} A uniformly small continuous
perturbation of the time-dependent coefficients
in the left-hand side of \eqref{basichyp} does not destroy its
unique solvability property. Hence, equation
\begin{equation}\label{aphyp}
\frac{d}{dt} v(t)- A(t)v + \<b(t)\cdot v\> f(y(t)) = h(t)
\end{equation}
has a unique uniformly bounded solution $v(t)$ given any uniformly bounded function $h(t)$;
moreover, for the corresponding operator $L: h\mapsto v$ estimate \eqref{clnrm} holds
(we assume that the constant $C_\Lambda$ in \eqref{clnrm} was taken with a margin of safety,
so all our small perturbations of the equation do not change $C_\Lambda$).
Note also, that given any function $\phi(t)$ that satisfies \eqref{fnu}, if we introduce a new time
$\tau=\phi(t)$ in the equation \eqref{aphphi} and a new function $v_{new}$ by the rule
$v_{new}(\phi(t))\equiv v(t)$, then the left-hand side of equation will be $O(\nu)$-close
to the left-hand side of \eqref{aphyp}. For sufficiently small $\nu$ this gives us the unique
solvability of \eqref{aphphi} and estimate \eqref{clnre} with $\alpha=0$.

Next, we note that a multiplication of the functions $v$ and $h$ in \eqref{aphphi} to any smooth
function of $t$ with uniformly small derivative just results in a uniformly small correction
to $A(\phi(t))$. This immediately shows the unique solvability of equation \eqref{aphphi}
in any weighted space with a sufficiently slowly growing weight; e.g. we obtain
\eqref{clnre} for all small $\alpha$.

In order to show the Lipshitz property of $L_\phi$ with respect to $\phi$, we note that
$$
v_2-v_1= L_{\phi_2}\left\{
(A(\phi_2)-A(\phi_1)) v_1-[\<b(\phi_2)\cdot v_1\> f(y(\phi_2))-\<b(\phi_1)\cdot v_1\> f(y(\phi_1))]\right\}.
$$
Now, since $A(\phi)$, $b(\phi)$, $f(y(\phi))$ are smooth - hence, Lipshitz - with respect to $\phi$,
and since $v_1(t)$ is uniformly bounded by \eqref{clnrm}, we immediately get \eqref{aplip} from \eqref{clnre}.
\end{proof}

Further we will use
\begin{equation}\label{bmu}
b(t)=\int_{-\infty}^{+\infty} y^*(t+s\mu) \xi(s) ds, \qquad
A(t)=\int_{-\infty}^{+\infty} f'(y(t+s\mu)) \xi(s) ds,
\end{equation}
where $\mu$ is a small constant and $\xi\geq 0$ is such that $\int_{-\infty}^{+\infty}\xi(s)ds=1$.
At $\mu=0$ we have $b\equiv y^*$ and $A\equiv f'(y)$; at small $\mu$ the functions $A(t)$ and $b(t)$
are close, respectively, to $f'(y(t))$ in $C^{r-1}$ and to $y^*(t)$ in $C^r$ (we have
$y^*(t)\in C^r$ as it satisfies equation \eqref{linvc}). Thus, uniformly for all $t$, we have
\begin{equation}\label{bfdot}
\begin{array}{l}
\<b(t)\cdot f(y(t))\> - 1:= c(t)=O(\mu),\\ \\
b'(t)+f'(y(t))^\top b(t)=o(1)_{\mu\to0},\qquad A(t)-f'(y(t))=o(1)_{\mu\to0}
\end{array}
\end{equation}
(see \eqref{dotdot},\eqref{linvc}). By taking $\xi\in C^\infty$ and such that
$\int_{-\infty}^{+\infty}|\xi'(s)|ds <\infty$, we will make $A(t)$ and $b(t)$ at $\mu\neq 0$
more smooth then $f'(y(t))$ and, respectively, $y^*(t)$, namely we will use $A$ which is at least $C^r$
and $b$ which is at least $C^{r+1}$; the price is that the last derivatives do not stay bounded
as $\mu\rightarrow0$, however we have en estimate
\begin{equation}\label{adsf}
\|A(t)\|_{C^r}=O(\mu^{-1}), \qquad \|b(t)\|_{C^{r+1}}=O(\mu^{-1}).
\end{equation}

The next proposition describes the way we coordinatize a small neighborhood of the
curve $w_k^0: y=y_k^0(t)$ in the $Y$-space ($y_k^0(t)$ is an orbit from the hyperbolic set
$\Lambda_k$ of the $k$-th subsystem of the uncoupled LDS \eqref{2.7}). Let $w_k: y=y_k(t)$
be a curve, $\gamma$-close to $w_k^0$ on some, finite or infinite, interval $I$ of $t$,
i.e. there exists a smooth time-reparametrization $\psi(t)$
such that $\|y_k^0(\psi(t))-y_k(t)\|_Y< \gamma$ at $t\in I$.

\begin{lemma}\label{stupid}
There exists $\bar\gamma>0$ (independent of the choice
of the orbit $y_k^0\in\Lambda_k$ and independent of k) such that if $\gamma<\bar\gamma$, then
there exists a uniquely defined on $I$ function $\phi(t)$ such that
$\phi=\psi+O(\gamma)$,  $\frac{d\phi}{d\psi}=1+o(1)_{\gamma\to0}$, and
\begin{equation}\label{vfix}
\<b_k(\phi)\cdot (y_k(t)-y_k^0(\phi))\>\equiv 0,
\end{equation}
where $b_k$ is given by \eqref{bmu} at some small $\mu$.
\end{lemma}
\begin{proof}
The derivative of the left-hand side of \eqref{vfix} with respect to $\phi$ at constant $y_k$ is
$\displaystyle \<b_k^\prime(\phi)\cdot (y_k-y^0_k(\phi))\>-
\<b_k(\phi)\cdot \dot y_k^0(\phi)\>=O(y_k-y^0_k(\phi))-\<b_k(\phi)\cdot f_k(y_k^0(\phi))\>$.
By \eqref{bfdot}, it is bounded away from zero, provided $y_k-y^0_k(\phi)$ is sufficiently small.
Thus, by the implicit function theorem, for any point $y_k$ from the (sufficiently small)
$\gamma$-neighborhood of the point $y_k^0(\varphi)$, we have a uniquely defined $\phi(y_k)$ which
satisfies \eqref{vfix} and condition $\phi(y_k^0(\psi))=\psi$. Moreover, $\phi$ depends smoothly on $y$
and the derivatives are uniformly bounded. So, as $\|y_k^0(\psi(t))-y_k(t)\|_Y< \gamma$, we also have
$\|\psi(t)-\phi(t)\|_Y = O(\gamma)$, as required (we denote $\phi(t):=\phi(y_k(t))$).
\end{proof}

Let us now proceed to the proof of the theorem. Let
$\yy^0(t)$ be an orbit from $\Lambda^\infty$. The sought invariant manifold
${\Bbb W}_{\yy^0,\varepsilon}$ consists of all solutions of the LDS \eqref{2.10}
which stay for all times in a small neighborhood of the manifold ${\Bbb W}_{\yy^0}^0$.
This means that, for every $k\in\Bbb Z$, the $k$-th component of $\yy(t)$ stays uniformly
close to the corresponding curve $w_k^0: y=y_k^0(\phi)$ in the $Y$-space. In other words
every trajectory $(\yy(t),\pp(t))\in{\Bbb W}_{\yy^0,\varepsilon}$ satisfies
\begin{equation}\label{yclose}
y_k(t)=y_k^0(\phi_k (t))+v_k(t), \ \ k\in\Bbb Z,
\end{equation}
where the functions $v_k(t)$ are uniformly small. By Lemma \ref{stupid}, we may always assume
that the parametrization $\phi_k(t)$ is chosen so that \eqref{vfix} is fulfilled.
By differentiating \eqref{vfix} with respect to $t$ we get
\begin{equation}\label{vfixd}
\<b_k(\phi_k(t))\cdot v_k'(t)\>\equiv -\phi_k'(t) \<b_k'(\phi_k(t))\cdot v_k(t)\>.
\end{equation}

Now, plugging \eqref{yclose} into the first equation of \eqref{2.10} gives
\begin{equation}\label{2.15}
v_k'(t)+\phi'_k(t)f_k(z_k(t))=f_k(z_k(t)+v_k(t))+\eb F_{\eb,k}(\zz+\vv(t),\pp(t)),
\end{equation}
where we denote $z_k(t):=y_k^0(\phi_k(t))$.
By multiplying both sides of this equation to $b(\phi_k(t))$, and taking \eqref{vfix},\eqref{vfixd} and
\eqref{bfdot} into account, we obtain the following equation for the evolution of $\phi_k$:
\begin{equation}\label{2.15p}
\phi_k'(t)=1+q_k(\vv,\Phi,\pp),
\end{equation}
where
\begin{equation}\label{qphik}
q_k:=\frac{\<[f_k(z_k+v_k)-f_k(z_k)]\cdot b_k(\phi_k)\>+ \<b'_k(\phi_k)\cdot v_k\>+\eb
\<F_{\eb,k}(\zz+\vv,\pp)\cdot b_k(\phi_k)\>}{1+c_k(\phi_k)-\<b_k'(\phi_k),v_k\>}.
\end{equation}
Equation for the $v$-components can now be obtained by plugging \eqref{2.15p}
into \eqref{2.15}:
\begin{equation}\label{2.15v}
v_k'(t)-A_k(\phi_k) v_k + \<b_k(\phi_k)\cdot v_k\>f_k(z_k)=Q_k(\vv,\Phi,\pp) - q_k(\vv,\Phi,\pp) f_k(z_k),
\end{equation}
where
\begin{equation}\label{rform}
Q_k:=f_k(z_k+v_k)-f_k(z_k)-A_k(\phi_k)v_k+\eb F_{\eb,k}(\zz+\vv,\pp).
\end{equation}
Equation for the evolution of $\pp(t)$ is given by the second equation of \eqref{2.10}:
\begin{equation}\label{2.16}
p_k'(t)=g_k(z_k+v_k)+\eb G_{\eb,k}(\zz+\vv,\pp).
\end{equation}

We remark that if we choose $\mu=0$ in \eqref{bmu}, equations \eqref{qphik} and \eqref{rform}
are simplified and reduce to
\begin{equation}\label{qm0}
q_k=\frac{\<Q_k\cdot b_k(\phi_k)\>}{1+\<y_k^*(\phi_k),f'_k(z_k)v_k\>},\quad
Q_k=f_k(z_k+v_k)-f_k(z_k)-f'_k(z_k)v_k+\eb F_{\eb,k}(\zz+\vv,\pp)
\end{equation}
(see \eqref{bfdot}). However, the functions $q_k$ and $Q_k$ will be only $C^{r-1}$ with respect
to $\Phi$, therefore we use small non-zero $\mu$ -- in order not to lose the last derivative
(and to be able to treat the case $r=1$).

By multiplying both sides of \eqref{2.15v} to $b_k(\phi_k(t))$
and using \eqref{2.15p},\eqref{bfdot} we find that
$$\frac{d}{dt} \<b_k(\phi_k(t))\cdot v_k(t)\>+ \<b_k(\phi_k(t))\cdot v_k(t)\>(1+O(\mu))=0.$$
This equation has only one bounded solution: $\<b_k(\phi_k(t))\cdot v_k(t)\>\equiv 0$; therefore,
since $b$ is uniformly bounded, we find that every uniformly bounded solution $\vv(t)$ of
the system \eqref{2.15p},\eqref{2.15v},\eqref{2.16} (with $k$ running all integer values)
satisfies \eqref{vfix}. Hence, it satisfies \eqref{2.15}. Thus, the solutions of system
\eqref{2.15p},\eqref{2.15v},\eqref{2.16} whose $\vv(t)$-component is uniformly small
give us all the solutions of system \eqref{2.10} which stay uniformly close
to the manifold ${\Bbb W}_{\yy^0}^0$ (i.e. all the solutions which comprise the sought
invariant manifold ${\Bbb W}_{\yy^0,\eb}$). We show below that for all small $\delta>0$
the solution of \eqref{2.15p},\eqref{2.15v},\eqref{2.16} for which
\begin{equation}\label{vsmal}
\|v_k(t)\|_Y\le \delta \qquad (k\in\Bbb Z, \; t\in\Bbb R)
\end{equation}
exists and is defined uniquely for any given initial condition $\Phi(0)$ and $\pp(0)$.

In order to prove the existence and uniqueness of the (small $\vv$) solution, we will show that
it can be obtained as a fixed point of a contracting operator on an appropriate space.
Namely, we consider the set $\Cal X_{\delta,\nu}$ of all functions $(\vv(t),\Phi(t),\pp(t))$
belonging to space $C_{loc}(\R,\Bbb Y\times\Psi\times\mathbb P)$
 such that \eqref{vsmal} and \eqref{fnu} hold
for all $k$ and $t$ for which the following norm is finite:
\begin{equation}\label{xnrm}
\|\vv,\Phi,\pp\|_\alpha=\sup_{k\in\Bbb Z, t\in\Bbb R} e^{-\alpha|t|}
\max\{\|v_k(t)\|, |\phi_k(t)|, \kappa\|p_k(t)\| \},
\end{equation}
where $\alpha>0$, and $\kappa>0$ is assumed to be sufficiently small. Obviously,
the set $\Cal X_{\delta,\nu}$ is a complete metric space with respect to that norm.

Note that in the limit $\lim_{\eb\to0,\vv\to0}$ the functions $q_k$, $Q_k$
given by \eqref{qphik},\eqref{rform} tend uniformly to zero for all $k\in\Bbb Z$, and in the limit
$\lim_{\mu\to+0}\lim_{\eb\to0,\vv\to0}$ their
first derivatives with respect to $\vv$, $\Phi$ and $\pp$ tend uniformly to zero too
(see \eqref{fgcr},\eqref{bfdot},\eqref{adsf}; the order of the limits is important: $\eb$ and $\vv$ first,
then $\mu$). The first derivative of the right-hand side of \eqref{2.16} with respect to $\pp$
is also uniformly small. Thus, if we rewrite system \eqref{2.15p},\eqref{2.15v},\eqref{2.16} as
\begin{equation}\label{2.22}
\left\{\begin{array}{ll}\displaystyle
v_k=L_{\phi_k}\left[Q_k(\vv,\Phi,\pp) - q_k(\vv,\Phi,\pp) f_k(z_k)\right],\\
\\ \displaystyle
\phi_k=\phi_k^0+t+\int_0^t q_k(\vv,\Phi,\pp) dt,\\
&\quad (k\in\Bbb Z)\\ \displaystyle
p_k=p_k^0+\int_0^t g_k(z_k+v_k) dt+\eb\int_0^t G_{\eb,k}(\zz+\vv,\pp) dt,
\end{array}\right.
\end{equation}
where $(\Phi^0,\pp^0)\in \mathbb P\times\Psi$ is arbitrary and the operator $L$ is defined by equation \eqref{aphphi},
then it is easy to check
that the right-hand side of \eqref{2.22}  (for every fixed $\Phi^0$ and $\pp^0$) defines a contracting operator
$\Cal T: \Cal X_{\delta,\nu}\to \Cal X_{\delta,\nu}$ for every exponential norm with a sufficiently small
weight $\alpha_0$. Namely, we first fix small $\nu$ and $\alpha_0$ such that the operators $L_{\phi_k}$
will all be defined and Lipshitz with respect to $\phi$ (see \eqref{aplip}, actually, we may fix $\nu$ of order $\eb$);
the operator of integration $\int_0^t (\cdot) dt$ is also Lipshitz in the $\alpha_0$-norm, with the Lipshitz
constant $\frac{1}{\alpha_0}$; then we choose a sufficiently small $\mu$ for which the Lipshitz
constants of $q_k$ and $Q_k$ can become small enough as $\eb$ and $\vv$ tend to zero; then we see
that one may choose $\kappa$ sufficiently small such that for all sufficiently small $\eb$ and $\delta$
the Lipshitz constant of the right-hand side of \eqref{2.22} on the space $\Cal X_{\delta,\nu}$ is
less than $1$, which means the operator $\Cal T$ is contracting indeed (we need to introduce the small
factor $\kappa$
in the definition of norm on $\Cal X_{\delta,\nu}$ because the derivative of $g_k$ with respect
to $v_k$ and $\phi_k$, though bounded, is not necessarily small). As at $\eb=0$ and $\vv=0$
the $\vv$-component of the image by $\Cal T$ vanishes, the contractivity of $\Cal T$ implies that
given any small $\delta$ the condition \eqref{vsmal} is invariant with respect to $\Cal T$
for all sufficiently small $\eb$; i.e., the $\Cal T \Cal X_{\delta,\nu}\subset \Cal X_{\delta,\nu}$.

By the Banach principle, the iterations by $\Cal T$ of any initial element from
$\Cal X_{\delta,\nu}$ converge to a uniquely defined limit in $\Cal X_{\delta,\nu}$, the fixed point of $\Cal T$.
Thus, we have shown that every solution which stays sufficiently close to the manifold
${\Bbb W}_{\yy^0}^0$ for all times can be found as the uniquely defined solution of \eqref{2.22}.
Therefore, the union of all such solutions comprises the sought invariant manifold ${\Bbb W}_{\yy^0,\eb}$
given by \eqref{wman} where the function $\Bbb V_k$ is defined by the  map
$(\Phi(0),\pp(0))\mapsto v_k(0)$. Note that this map (hence the manifold ${\Bbb W}_{\yy^0,\eb}$) is
Lipshitz continuous, since the contracting operator $\Cal T$ is Lipshitz continuous with respect
to $(\Phi^0,\pp^0)$. We omit the proof of the smoothness of this map, as it is completely standard (yet
laborious): one may show that the operator $\Cal T$ is smooth on a scale of Banach spaces corresponding to different
weighted $\alpha$-norms on $\Cal X_{\delta,\nu}$ (cf. \cite{481,41}), or alternatively check, by fiber-contraction
arguments, that the iterations of by $\Cal T$ of an initial element from
$\Cal X_{\delta,\nu}$ converge to the fixed point of $\Cal T$ uniformly along with the derivatives with respect
to $(\Phi^0,\pp^0)$ (cf. \cite{30,3001}).

In order to finish the proof of the theorem, it remains to show estimate \eqref{2.25}.
The $C^r$-part is obvious, as $\vv=0$ solves \eqref{2.22} at $\varepsilon=0$, and the fixed point
of a contracting operator which depends on a parameter continuously must depend on the same parameter
continuously. To show the $C^{r-1}$-estimate, we note that when the right-hand side of \eqref{2.22} depends
smoothly on some parameter, the solution must also be smooth with respect to the same parameter.
In particular, if we rewrite system \eqref{2.10} as
$$\yy'(t)=\ff(\yy)+\sigma\FF_\eb(\yy,\pp),\qquad \pp'(t)=\ggg(\yy)+\sigma\GG_\eb(\yy,\pp),$$
then $\vv(t)$ will depend $C^r$-smoothly on $\sigma$ as well, which immediately gives
\eqref{2.25} if we note that $\vv(t)=0$ at $\sigma=0$ and plug $\sigma=\eb$ back.
\end{proof}

\begin{remark}\label{rem0k}
The theorem also remains true if the number of systems coupled in the LDS is finite, i.e. if
the index $k$ runs a finite set instead of $\Bbb Z$. Then the range of $\eb$ values
for which the corresponding invariant manifolds exist will be independent on the number $N$
of systems in the LDS -- provided the constant $C$ in the bound \eqref{fgcr} on the norm of the
coupling terms is independent of $N$. Note that condition \eqref{fgcr} does not requires
that the coupling is local, it just means that the ``total coupling strength'' for each
subsystem in the LDS is bounded independently of the total number $N$ of subsystems involved.
\end{remark}

\begin{remark}\label{rem1k}
Given any symmetry in system \eqref{2.10}, if the set $\Lambda^\infty$ obeys the same symmetry,
then system of invariant manifolds $\Bbb W_{\yy^0,\eb}$ inherits the symmetry for all small $\eb$
-- by uniqueness.
A basic example of such symmetry is invariance with respect to {\em spatial translation} $k\to k+1$ (in
this case, the coupling terms $\FF,\GG$ are shift-invariant,
the individual ODE's \eqref{2.7} are the same
for all $k$, and the sets $\Lambda_k$ should be chosen the same).
\end{remark}

Theorem \ref{Th2.1} allows us to construct a huge number of
special solutions of the weakly coupled LDS \eqref{2.10}. Indeed,
for every $\yy^0\in\Lambda^\infty$, one can construct the
associated manifold $\Bbb W_{\yy^0,\eb}$ and then, for every
$(\Phi_0,\pp_0)\in\Psi\times\Bbb P$, there exists a solution
$(\yy(t),\pp(t))$ in the form
\begin{equation}\label{2.31}
\yy(t)=\yy^0(\Phi(t))+\Bbb V_{\yy^0}(\Phi(t),\pp(t))
\end{equation}
where the functions $(\Phi(t),\pp(t))$ solve the reduced problem
on the center manifold (see \eqref{2.15p},\eqref{2.16}):
\begin{equation}\label{2.32}
\begin{cases}
\Phi'(t)=1+\qq(\Bbb V_{\yy^0}(\Phi,\pp),\Phi,\pp),\\
\pp'=\ggg(\yy^0(\phi)+\Bbb V_{\yy^0}(\Phi,\pp))+\eb\GG_\eb(\yy^0(\Phi)+\Bbb V_{\yy^0}(\Phi,\pp),\pp),\\
\Phi(0)=\Phi_0,\ \ \pp(0)=\pp_0.
\end{cases}
\end{equation}
It is interesting to have an expansion in powers of $\eb$ for the system on the invariant manifold.
In order to do this we need a sufficient smoothness of the right-hand sides: for instance, to find
the first order in $\eb$ approximation to \eqref{2.32} we assume the original system to be at least $C^2$
with respect to all variables and $\eb$. In this case we may take $\mu=0$ in formulas \eqref{bmu},
so the function $\qq$ will be given by \eqref{qm0}. As $\Bbb V=O(\eb)$ by \eqref{2.25}, we immediately
obtain the first-order in $\eb$ approximation to the $\phi$-equation:
\begin{equation}\label{ph1e}
\phi_k'(t)=1+\eb \<F_{k,0}(\yy^0(\Phi),\pp)\cdot y_k^*(\phi_k)\>,
\end{equation}
where $y_k^*(s)$ is the uniquely defined bounded solution of
$$\frac{d}{ds} y_k^*(s) = -f_k'(y_k^*(s))^\top y_k^*(s), \qquad \<y^*_k(s)\cdot f_k(y_k^0(s))\>\equiv1.$$
Formula \eqref{ph1e} describes the evolution of phases on the invariant manifolds $\Bbb W$ and can be useful
in the study of phase synchronization in coupled chaotic systems (see e.g. \cite{Pik}).

To obtain the approximate $\pp$-equation, we need the first-order approximation to $\vv$. By expanding
the first equation in \eqref{2.22} in $\eb$, we find that
$$v_k=\eb u_k(0)+o(\eb),$$
where the function $u_k(t)$ is given by
$$u_k=L_{\phi_k(t)}
\left[F_{k,0}(\yy^0(\Phi(t)),\pp(t))-
\<F_{k,0}(\yy^0(\Phi(t)),\pp(t))\cdot y_k^*(\phi_k(t))\> f_k(y_k^0(\phi_k(t)))\right].$$
Since, by \eqref{2.32}, $\phi_k(t)$ is for all $k$ uniformly $O(\eb)$-close to
$\phi_k(0)+t$ in the exponential
$\alpha$-norm with $\alpha>0$, and $p_k(t)$ is for all $k$ uniformly $O(\eb)$-close to
$p_k(0)+\int_0^tg_k(y_k^0(s+\phi_k(0)))ds$, also in the exponential
$\alpha$-norm, it follows from the Lipshitz property of the operator $L_\phi$ (see
\eqref{aplip},\eqref{clnre}) that
$$v_k=\eb w_k(\Phi,\pp)+o(\eb),$$
where, given any constant $\Phi$ and $\pp$, we denote as $w_k(\Phi,\pp)$ the value at $t=0$
of the uniquely defined bounded solution $w(t)$ of the equation
$$
\begin{array}{l}
\frac{d}{dt} w(t)- f'_k(y_k^0(t+\phi_k))w + \<y_k^*(t+\phi_k)\cdot v\> f_k(y_k(t+\phi_k))=\\
\qquad =\;\;\Cal F_k - \<\Cal F_k\cdot y_k^*(t+\phi_k)\> f_k(y_k^0(t+\phi_k))),\quad \mbox{where}\\
\Cal F_k:=F_{k,0}(\yy^0(t+\Phi),\pp+\int_0^t\ggg(\yy^0(s+\Phi))ds).
\end{array}$$
By plugging the above formula for $v_k$ into the $\pp$-equation of \eqref{2.32} and dropping all
$o(\eb)$-terms, we find that the first-order approximation to the $\pp$-equation is
\begin{equation}\label{pp1e}
p_k'(t)=g_k(\yy^0_k(\phi_k)+\eb w_k(\Phi,\pp))+\eb G_{k,0}(\yy^0(\Phi),\pp).
\end{equation}

\begin{remark}\label{rem2k} Let the uniformly-hyperbolic sets $\Lambda_k$ be compact and
{\em locally-maximal}, i.e. there exists $\gamma^0>0$ (independent of $k$) such that for each $k$
every orbit of \eqref{2.7}, which stays in the $\gamma^0$-neighborhood of $\Lambda_k$ for all $t$,
belongs to $\Lambda_k$ itself. Then, for $\eb$ sufficiently small,
every solution of the coupled LDS \eqref{2.10}
whose $\yy$-component stays for all times in a small neighborhood
of $\Lambda^\infty$ belongs to one of the manifolds $\Bbb W_{\yy^0,\eb}$. Indeed, given any $k$
the $k$-th component $y_k(t)$ of such solution must be close, after some reparametrization of time, to
a $\gamma${\em-orbit} $\;\tilde y(\varphi(t))$, which is a countable union of consecutive pieces
$\tilde y(\varphi)|_{\phi\in[\phi_j,\phi_{j+1})}$ of orbits from the set $\Lambda_k$ such that
$\|\tilde y(\varphi_j)-\tilde y(\varphi_j-0)\|\leq\gamma$, for some small $\gamma$. It is known that
when $\Lambda_k$ is locally-maximal, any $\gamma$-orbit is shadowed by a true orbit, i.e. there exists
an orbit in $\Lambda_k$ which is $O(\gamma)$-close to $\tilde y(\varphi)$ (after a reparametrization
of time). Thus, our solution $\yy(t)$ of the coupled LDS stays for all times close to a (time-reparametrized)
orbit $\yy^0\in\Lambda^\infty$, i.e. we can write it in the form \eqref{yclose}, and we showed in Theorem
\ref{Th2.1} that every such solution belongs to the invariant manifold $\Bbb W_{\yy^0,\eb}$.
\end{remark}

\section{A theorem on asymptotic phase}\label{s2.5}
In this Section we compare the behavior of orbits of the LDS \eqref{2.10} which belong to
different invariant manifolds $\Bbb W_{\yy^0,\eb}$. We start with the analysis of the dependence of the
invariant manifold $\Bbb W_{\yy^0,\eb}$ on the choice of the trajectory $\yy^0\in\Lambda^\infty$.
Clearly, $\Bbb W_{\yy^0,\eb}$ depends on $\yy^0$ continuously: namely, the
function $\Bbb V_{\yy^0,\eb}(\Phi,\pp)$ (hence - the function $\Bbb U_{\yy^0,\eb}(\Phi,\pp)$)
in \eqref{wman} is found via an application of the contraction mapping principle, and the corresponding
contracting operator (the operator $\Cal T$ defined by the right-hand side of \eqref{2.22})
depends continuously on $\yy^0$ in some exponential weighted norm, so on any
bounded set of values of $\Phi$ and $\pp$ the functions $\Bbb U_{\yy^{0,1},\eb}$ and
$\Bbb U_{\yy^{0,2},\eb}$ will be uniformly close provided the trajectories $\yy^{0,1}$ and $\yy^{0,2}$
are sufficiently close in the weighted norm. We need, however, a different statement about the
closeness of $\Bbb U_{\yy^{0,1},\eb}$ and $\Bbb U_{\yy^{0,2},\eb}$.
Note that though the manifold $\Bbb W_{\yy^0,\eb}$ is defined uniquely (as the set of all solutions
that for all $t$ stay uniformly close to the manifold $\Bbb W_{\yy^0}^0$ defined by \eqref{wman0}),
the function $\Bbb U_{\yy^0,\eb}$ in \eqref{wman} is defined up to an (arbitrary) reparametrization
of the space $\Psi$ of phases $\phi_k$. Therefore, when comparing functions $\Bbb U$ corresponding
to two different trajectories $\yy^0$ (as we do it below), we should describe how the corresponding
parametrization choices agree with each other.

In order to do so we recall the construction used
in the proof of Theorem \ref{Th2.1}. Take any orbit $\{\yy=\yy^0(t)\}\in\Lambda^\infty$;
its $k$-th component $y_k^0(t)$ defines a smooth curve $w^0_k$ in the space $Y$. Take any other
curve $w: \{y=y_k(t)\}$ in $Y$. By Lemma \ref{stupid},
there exists $\bar\gamma>0$ (independent of the choice of the
curves) such that if $w$ stays in the $\bar\gamma$-neighbourhood of $w_k^0$ for a certain
interval of time, then, for every $t$ from this interval, condition \eqref{vfix} defines
the projection $y_k^0(\phi_k(t))$ of the point $y_k(t)\in w$ onto the curve $w_k^0$ uniquely.
We will call $\phi_k(t)$ the
phase {\em relative to} $\yy^0$. If we have two orbits, $\yy^{0,1}$ and $\yy^{0,2}$, from
$\Lambda^\infty$, and these orbits are $\gamma$-close ($\gamma<\bar\gamma$) on some time interval,
then for every curve $y=y_k(t)$ which stays at the distance less than $\gamma$ from both $w_k^{0,1}$
and $w_k^{0,2}$ on this time interval we have two phases, $\phi_k^1(t)$ and $\phi_k^2(t)$,
relative to $\yy^{0,1}$ and $\yy^{0,2}$ respectively. By Lemma \ref{stupid} (with $\psi$ standing for $\phi_k^2$
and $\phi$ for $\phi_k^1$), these two phases
are related by a close to identity diffeomorphism: $\phi_k^2(t)=\eta_k(\phi_k^1(t))$
where $\eta_k^\prime(\phi)=1+o(1)_{\gamma\to0}$. For a solution in the invariant manifold
$\Bbb W_{\yy^0,\eb}$ the canonical phases $\phi_k$ (which we used before) are the phases relative
to $\yy^0$. However, if two orbits $\yy^{0,1}$ and $\yy^{0,2}$ are $\gamma$-close ($\gamma<\bar\gamma$)
on some time interval, then for solutions in, say, $\Bbb W_{\yy^{0,j},\eb}$ the phase $\varphi_k$
relative to $\yy^{0,1}$ is also defined on this interval, along with the canonical phase $\phi_k$
relative to $\yy^{0,2}$.

\begin{lemma}\label{Cor2.1} Let the assumptions of Theorem \ref{Th2.1} hold.
Then there exists $\bar\gamma>0$, $\alpha>0$ and $C>0$ such that, for all small $\eb$, given any $T>T_0>0$,
if any two orbits $\yy^{0,1}$ and $\yy^{0,2}$ from $\Lambda^\infty$ satisfy
\begin{equation}\label{blizk}
\sup_{t\in [-T,T]}\|\yy^{0,1}(t)-\yy^{0,2}(t)\|_{\Bbb Y}<\gamma,
\end{equation}
for some $\gamma<\bar\gamma$, then there exists a uniformly
close to identity diffeomorphism $\eta$ such that
\begin{equation}\label{2.34}
\|\Bbb U_{\yy^{0,1}}(\Phi_0,\pp_0)-\Bbb U_{\yy^{0,2}}(\eta(\Phi_0),\pp_0)\|_{\Bbb Y}\le C e^{-\alpha (T-T_0)}
\end{equation}
for all $\pp_0\in\Bbb P$ and $\Phi_0$ such that
\begin{equation}\label{bndp}
\|\Phi_0\|_\Psi\leq T_0.
\end{equation}
\end{lemma}

\begin{proof} Let $(\yy^1(t),\pp^1(t))$ and $(\yy^2(t),\pp^2(t))$ be the orbits on the invariant manifolds,
respectively,
$\Bbb W_{\yy^{0,1}}$ and $\Bbb W_{\yy^{0,2}}$ such that $(\yy^1(t),\pp^1(t))$ corresponds to the initial condition
$\Phi(0)=\Phi_0$ and $\pp(0)=\pp_0$, and $(\yy^2(t),\pp^2(t))$ corresponds to the initial condition
$\Phi(0)=\eta(\Phi_0)$ and $\pp(0)=\pp_0$, where $\eta_k:\varphi_k\mapsto\phi_k$ is the close
to identity diffeomorphism which sends the phases relative to $\yy^{0,1}$ to the phases relative to
$\yy^{0,2}$.

Let $\phi_k(t)$ be canonical phases of $\yy^1(t)$ and let $\vv(t):=\yy^1(t)-\yy^{0,1}(\Phi(t))$,
so
\begin{equation}\label{uvv1}
\vv(0)=\Bbb U_{\yy^{0,1}}(\Phi_0,\pp_0)-\yy^{0,1}(\Phi_0)
\end{equation}
(see \eqref{2.31}). Let $\varphi_k(t)$ denote the phase
of $y_k^2(t)$ relative to $\yy^{0,1}$, and let
$\uu(t):=\{u_k(t)\}_{k\in\Bbb Z}$, where $u_k(t)=y_k^2(t)-y_k^{0,1}(\varphi_k(t))$. By construction,
$\varphi_k(0)=\phi_k(0)$ for all $k$, so
\begin{equation}\label{uvv2}
\uu(0)=\Bbb U_{\yy^{0,2}}(\eta(\Phi_0),\pp_0)-\yy^{0,1}(\Phi_0).
\end{equation}

Denote $x_k^1(t):=(v_k(t),\phi_k(t),p_k^1(t))$ and $x_k^2(t):=(u_k(t),\varphi_k(t),p_k^2(t))$;
as we showed in the proof of Theorem \eqref{Th2.1} the functions $x_k^j(t)$ satisfy the same system
\eqref{2.15p},\eqref{2.15v},\eqref{2.16}: $x_k^1(t)$ satisfies this system for all $t$, while
$x_k^2(t)$ satisfies it for all $t$ for which $\uu(t)$ remains small. As the orbit $\yy^2(t)$ belongs
to the invariant manifold $\Bbb W_{\yy^{0,2},\eb}$, it stays close to $\yy^{2,0}(\Phi^2(t))$ for
all times (where $\Phi^2$ is the canonical phase of $\yy^2$),
so by \eqref{blizk} the distance $\|\uu(t)\|$ between $\yy^2(t)$ and its projection to
$\Bbb W_{\yy^{0,1}}^0$ will remain small for all times such that $\|\Phi^2(t)\|\leq T$. Since
the time derivative of $\Phi$ is bounded (see \eqref{2.31}), we have from\eqref{bndp} the required
smallness of $\uu(t)$ for all $|t|\leq S$ where
\begin{equation}\label{stt0}
S=O(T-T_0+1).
\end{equation}

Outside this time-interval we cannot guarantee that the phases $\varphi_k(t)$ are well-defined, therefore we
modify $x_k^2(t)$ at $|t|\geq S-1$. Namely, we consider the functions
$x_k^3(t)=(v_k^3(t),\phi_k^3(t),p_k^3(t))$ defined by the following rule:
\begin{equation}\label{ok}
\begin{array}{l}
v_k^3(t)=\theta_0(t))u_k(t), \qquad p_k^3(t)=\theta_0(t) p_k^2(t),\\ \\
\phi_k^3(t)=\theta_0(t)\varphi_k(t)\;+\quad
\theta_-(t)[\varphi_k(-S+1)+\varphi_k'(-S+1)(t+S-1)]+\\\qquad\qquad
\;+\;\;\theta_+(t)[\varphi_k(S-1)+\varphi_k'(S-1)(t-S+1)],
\end{array}
\end{equation}
where $\theta_\pm(t)$ are smooth functions $\Bbb R^1\to[0,1]$ such that $\theta_-(t)$ equals to $1$ at
$t\leq -S$ and to $0$ at $t\geq -S+1$, while $\theta_+(t)$ equals to $1$ at
$t\geq S$ and to $0$ at $t\leq S-1$, and $\theta_0:=1-\theta_+-\theta_-$.
Note that it follows from \eqref{ok} that $v_k^3(t)$ is
uniformly small for all $t\in\Bbb R$ since $\uu(t)$ is uniformly small for all $|t|\leq S$.

Since both $\xx^1(t)$ and $\xx^2(t)$ satisfy system \eqref{2.15p},\eqref{2.15v},\eqref{2.16}
at $t\in[t_1+S,t_2-S]$, the function $\xx^3(t)$ satisfies the same system
with a uniformly bounded correction to the right-hand sides which is localized at
$|t|\in[S-1,S]$ (and which is denoted below as $\rho$). Since the initial conditions
in $\Phi$ and $\pp$ coincide for $\xx^1$ and $\xx^2$ by construction (recall that we choose $\yy^1(t)$
and $\yy^2(t)$ such that $\varphi_k(0)=\phi_k(0)$), we find that $\xx^3(t)_{t\in[-\infty,\infty]}$
satisfy the following equation (a perturbation of \eqref{2.22})
\begin{equation}\label{2.22j}
\left\{\begin{array}{ll}\displaystyle
v_k^3=L_{\phi_k^3}\left[Q_k(\vv^3,\Phi^3,\pp^3) - q_k(\vv^3,\Phi^j,\pp^3) f_k(z_k^3)+\rho_{k1}\right],\\
\\ \displaystyle
\phi_k^3=\phi_k(0)+t+\int_0^t \left[q_k(\vv^3,\Phi^3,\pp^3)+\rho_{k2}\right] dt,\\
&\quad (k\in\Bbb Z)\\ \displaystyle
p_k^3=p_k(0)+\int_0^t \left[g_k(z_k^3+v_k^3)+\rho_{k3}\right] dt
+\eb\int_0^t G_{\eb,k}(\zz^3+\vv^3,\pp^3) dt,
\end{array}\right.
\end{equation}
where $z_k^3:=y^{0,1}_k(\phi_k^3)$ and the perturbations $\rho_k(t)$ satisfy
$\|\rho\|_{\alpha_0}=O(e^{-\alpha_0 S})$.

Recall that $x_k^2(t)=(u_k(t),\varphi_k(t),p_k^2(t))$ satisfies system
\eqref{2.15p},\eqref{2.15v},\eqref{2.16} at $|t|\leq S$, and $\uu(t)$ is uniformly small on this interval
(provided $\gamma$ and $\eb$ are small enough). The smallness of $\uu$ and $\eb$ implies the smallness
of the functions $q_k$ in the right-hand side of the equation \eqref{2.15p} for the phases $\varphi_k$,
therefore $\sup_{|t|\leq S}|\varphi_k^\prime(t)-1|$ is uniformly small for all $k$. By \eqref{ok}, we find then
that $\sup_{t\in\Bbb R}|\phi_k^{3\prime}(t)-1|$ is also uniformly small. This guarantees
that the operator $L_{\phi_k^3}$ is defined and Lipshitz in the $\alpha_0$-norm
(see comments after \eqref{2.22} in the proof of Theorem \ref{Th2.1}). Since operator $L_\phi$
is Lipshitz in the exponential $\alpha_0$-norm, and so is the operator of integration
$\int_0^t (\cdot) dt$, we may rewrite \eqref{2.22j} as
$$\xx^3=\Cal T \xx^3 + O(e^{-\alpha_0 S})_{\alpha_0},$$
where $\Cal T$ is the operator defined by the right-hand side of \eqref{2.22}, i.e.
$\displaystyle \xx^1=\Cal T \xx^1$,
and we immediately get that
$$\|\xx^3-\xx^1\|_{\alpha_0}=O(e^{-\alpha_0 S}),$$
since the operator $\Cal T$ is contracting (in the norm given by \eqref{xnrm}; note that, as we have shown
in the proof of Theorem \eqref{Th2.1}, in order
to have contraction, both $\vv^1$ and $\vv^3$ must be uniformly small, i.e. must satisfy \eqref{vsmal}
with a sufficiently small $\delta$, and this property indeed holds true when $\eb$ and $\gamma$ are
sufficiently small).

In particular (since $\vv^3(0)=\uu(0)$), we have $\displaystyle \|\vv(0)-\uu(0)\|_{\Bbb Y}=O(e^{-\alpha_0 S})$,
and the lemma follows from \eqref{uvv1},\eqref{uvv2},\eqref{stt0}.
\end{proof}
\begin{remark}\label{remshift}
By a shift of time, we obtain that if
\begin{equation}\label{blizks}
\sup_{t\in [T_1,T_2]}\|\yy^{0,1}(t)-\yy^{0,2}(t)\|_{\Bbb Y}<\gamma,
\end{equation}
then
\begin{equation}\label{2.34s}
\|\Bbb U_{\yy^{0,1}}(\Phi_0,\pp_0)-\Bbb U_{\yy^{0,2}}(\eta(\Phi_0),\pp_0)\|_{\Bbb Y}\le C e^{-\alpha T)}
\end{equation}
for all $\pp_0\in\Bbb P$ and $\Phi_0$ such that for all $k$ the components $\phi_k$ of $\Phi_0$ satisfy
\begin{equation}\label{bndps}
T_1+T\leq \phi_k\leq T_2-T.
\end{equation}
\end{remark}
\medskip

We may now prove the following theorem, crucial for the next Section.
\begin{theorem}\label{Prop2.trace} Let the assumptions of Theorem
\ref{Th2.1} hold. Then there exists $\alpha>0$ and $\bar\gamma>0$ such that for all sufficiently small
$\eb$ and all $\gamma\in(0,\bar\gamma)$,
for every two trajectories $\yy^{0,1}$ and $\yy^{0,2}$ from $\Lambda^\infty$ satisfying the condition
\begin{equation}\label{2.expcl}
\sup_{t\geq t_0} \|\yy^{0,1}(t)-\yy^{0,2}(t)\|_{\Bbb Y}\leq \gamma,
\end{equation}
given any solution $(\yy^1(t),\pp^1(t))$ from the invariant manifold
$\Bbb W_{\yy^{0,1},\eb}$, there exists a unique solution $(\yy^2(t),\pp^2(t))$ from the invariant manifold
$\Bbb W_{\yy^{0,2},\eb}$ such that
\begin{equation}\label{2.trace}
\|\yy^1(t)-\yy^2(t)\|_{\Bbb Y}+\|\pp^1(t)-\pp^2(t)\|_{\Bbb P}\le C(\gamma) e^{-\alpha (t-t_0)}, \quad t\geq t_0.
\end{equation}
The factor $C(\gamma)$ tends to zero as $\gamma\to0$.
\end{theorem}
\begin{remark}\label{plm}
Absolutely analogously,
for every two trajectories $\yy^{0,1}$ and $\yy^{0,2}$ from $\Lambda^\infty$ satisfying the condition
\begin{equation}\label{2.expclm}
\sup_{t\leq t_0} \|\yy^{0,1}(t)-\yy^{0,2}(t)\|_{\Bbb Y}\leq \gamma,
\end{equation}
given any solution $(\yy^1(t),\pp^1(t))$ from the manifold
$\Bbb W_{\yy^{0,1},\eb}$, there exists a unique solution $(\yy^2(t),\pp^2(t))$ from
$\Bbb W_{\yy^{0,2},\eb}$
such that
\begin{equation}\label{2.tracem}
\|\yy^1(t)-\yy^2(t)\|_{\Bbb Y}+\|\pp^1(t)-\pp^2(t)\|_{\Bbb P}\le C(\gamma) e^{-\alpha |t-t_0|}, \quad t\leq t_0.
\end{equation}
\end{remark}
\medskip
\begin{proof}
As we explained in Lemma \ref{Cor2.1}, condition \eqref{2.expcl}
(which is an analogue of condition \eqref{blizk} for the case of infinite
time interval) implies that for all sufficiently small $\eb$, for
any two solutions $(\yy^1(t),\pp^1(t))$ and $(\yy^2(t),\pp^2(t))$ from the invariant manifolds,
respectively, $\Bbb W_{\yy^{0,1},\eb}$ and $\Bbb W_{\yy^{0,2},\eb}$, for
all $t\geq t_0$ and every $k\in\Bbb Z$ we have well-defined projections of the points
$y^1_k(t)$ and $y_k^2(t)$ onto the curve $y=y_k^{0,1}(\varphi_k)$ in the space $Y$. The position
of the projection point is defined by its phase $\varphi_k$, so we have two phases (relative
to the same orbit $\yy^{0,1}$) defined for all $t\geq t_0$: $\varphi_k^1(t)$ for the point $y^1_k(t)$
and $\varphi^2_k(t)$ for $y_k^2(t)$. Thus,
\begin{equation}\label{ykvi}
y^i_k(t)=y^{0,1}_k(\varphi^i_k(t))+v^i_k(t), \;\mbox{ where }\;
\<b_k(\varphi^i_k(t))\cdot v^i_k(t)\>\equiv 0\quad t\geq t_0;
\end{equation}
here $b_k$ is given by \eqref{bmu} with $y^*=y^{*,1}_k$.

As the solution $(\yy^1(t),\pp^1(t))$ belongs to the invariant manifold
$\Bbb W_{\yy^{0,1},\eb}$ associated with the orbit $\yy^{0,1}$
relative to which the phase is defined, the phases
$\varphi_k^1$ are just the canonical phases $\phi_k^1$. For the solution $(\yy^2(t),\pp^2(t))$, as we explained
in the introduction to Lemma \ref{Cor2.1},
the phases $\varphi_k^2$ are related to the canonical phases $\phi_k^2$ by a close
to identity diffeomorphism $\eta_k:\varphi_k^2\mapsto \phi_k^2$ at $t\geq t_0$; so,
$\phi_k^2(t)-\varphi_k^2(t)$ is uniformly small for all $t\geq t_0$.

Formula \eqref{ykvi} is identical to \eqref{yclose},\eqref{vfix}, hence
(see the proof of Theorem \ref{Th2.1})
the functions $(\vv^i(t),\varphi^i(t),\pp^i(t))$, both for $i=1$ and $i=2$, solve the same
system \eqref{2.15p},\eqref{2.15v},\eqref{2.16} (where one should replace $\phi$ with $\varphi$
and $z_k$ with $y_k^{0,1}(\varphi_k)$), for all $t$ for which $\vv^i(t)$ remains small.
As the solution $(\yy^1(t),\pp^1(t))$ belongs
to the invariant manifold $\Bbb W_{\yy^{0,1},\eb}$, we have that $\yy^1(t)$ stays close to
$\yy^{1,0}(\Phi^1(t))$ for all times, which guarantees the smallness of $\vv^1$. The smallness of
$\vv^2(t)$ at all $t\geq t_0$ follows from the fact that $(\yy^2(t),\pp^2(t))$ belongs
to the invariant manifold $\Bbb W_{\yy^{0,2},\eb}$, hence $y^2_k(t)$ stays close to
$y^{2,0}_k(\phi^2_k(t))$ for all times, and because of the uniform closeness of $\phi^2_k(t)$ to
$\varphi^2_k(t)$ and $y_k^{2,0}$ to $y_k^{1,0}$ at $t\geq t_0$ we obtain the uniform closeness
of $y^2_k(t)$ to $y^{1,0}_k(\varphi^2_k(t))$ at $t\geq t_0$.

Thus, we have
\begin{equation}\label{szs}
\frac{d}{dt}\varphi_k^i=1+q_k(\vv^i,\varphi^i,\pp^i),\qquad
\frac{d}{dt}p_k^i= h_k(\vv^i,\varphi^i,\pp^i)\qquad\qquad (k \in\Bbb Z),
\end{equation}
where $q_k$ is given by \eqref{qphik} (the only important thing for us is that
$q_k$ is uniformly small along with its first derivatives), and
\begin{equation}\label{hkforno}
h_k(\vv,\varphi,\pp):=
g_k(y_k^{0,1}(\varphi_k)+v_k)+\eb G_{\eb,k}(\yy^{0,1}(\varphi)+\vv,\pp).
\end{equation}
By Theorem \eqref{Th2.1}, since the solutions $(\yy^i(t),\pp^i(t))$ belong
to the respective invariant manifolds $\Bbb W_{\yy^{0,i},\eb}$, we may put
\begin{equation}\label{vinvma}
\vv^i(t)=\Bbb V^i(\varphi^i(t),\pp^i(t)),
\end{equation}
in equations \eqref{szs}, where $\tilde{\Bbb V}^i$ are certain functions of $(\varphi,\pp)$
with the Lipshitz constant uniformly small. Namely,
the function $\tilde{\Bbb V}^1$ is just the function $\Bbb V^1$ that defines the manifold
$\Bbb W_{\yy^{0,1},\eb}$ by \eqref{wman}, while the function $\tilde{\Bbb V}^2$ is given by
$\tilde{\Bbb V}(\varphi,\pp)=\Bbb V^2(\eta(\varphi),\pp)+\yy^{0,2}(\eta(\varphi))-\yy^{0,1}(\varphi)$,
where $\eta$ is the diffeomorphism which sends the phase relative to $\yy^{0,1}$ to the phase
relative to $\yy^{0,2}$; the required Lipshitz property
of $\tilde{\Bbb V}^2$ follows from the Lipshitz property of $\Bbb V$ and $\eta$. Note that by \eqref{wman}
\begin{equation}\label{tlvtlu}
\tilde{\Bbb V}^1(\varphi,\pp)=\Bbb U^1(\varphi,\pp)-\yy^{0,1}(\varphi),\qquad
\tilde{\Bbb V}^2(\varphi,\pp)=\Bbb U^2(\eta(\varphi),\pp)-\yy^{0,1}(\varphi),
\end{equation}
hence, by Remark \ref{remshift}, when $\varphi_k\to+\infty$ uniformly
for all $k$, we have
\begin{equation}\label{1mky}
\tilde{\Bbb V}^2(\varphi,\pp)-\tilde{\Bbb V}^1(\varphi,\pp)=O(e^{-\alpha'\|\varphi\|_\Psi})
\end{equation}
for some $\alpha'>0$.

It follows from \eqref{1mky},\eqref{vinvma},\eqref{ykvi} that we will have the required
exponential decay of $\|\yy^1(t)-\yy^2(t)\|$ if the difference
between the corresponding solutions
$(\varphi^2(t),\pp^2(t))$ and $(\varphi^1(t),\pp^1(t))$ of \eqref{szs},\eqref{vinvma}
tends exponentially to zero as $t\to+\infty$. Given $(\varphi^1(t),\pp^1(t))$,
the sought, tending to it solution $(\varphi^2(t),\pp^2(t))$ corresponds to the fixed
point of the operator $(\phi(t),\pp(t))_{t\geq0}\mapsto(\bar\phi(t),\bar\pp(t))_{t\geq0}$
defined by the following equation:
\begin{equation}\label{2.cor}
\begin{aligned}
\bar\phi(t)=\int_t^\infty
\left[\qq(\tilde{\Bbb V}^1(\varphi^1,\pp^1),\varphi^1,\pp^1)-
\qq(\tilde{\Bbb V}^2(\varphi^1+\phi,\pp^1+\pp),\varphi^1+\phi,\pp^1+\pp)\right] dt,\\
\bar\pp(t)=\int_t^\infty\left[\hh(\tilde{\Bbb V}^1(\varphi^1,\pp^1),\varphi^1,\pp^1)-
\hh(\tilde{\Bbb V}^2(\varphi^1+\phi,\pp^1+\pp),\varphi^1+\phi,\pp^1+\pp)\right] dt,
\end{aligned}
\end{equation}
where we denote $\displaystyle \phi(t):=\varphi^2(t)-\varphi^1(t), \ \ \pp(t):=\pp^2(t)-\pp^1(t)$.
Thus, it remains to prove the existence and uniqueness
of the fixed point of this operator in the space of exponentially decreasing
functions, and also to show that this fixed point tends to zero as $\gamma\to0$.

In order to do this, we first note that
because $\varphi^1(t)$ and $\varphi^2(t)$ grow within linear bounds with time,
estimate \eqref{1mky} along with the boundedness of the Lipshitz constants
of the functions $\qq$ and $\hh$ implies that for some $\alpha>0$ the operator \eqref{2.cor}
takes exponentially decreasing functions
$(\phi(t),\pp(t))=O(e^{-\alpha t})$ into functions $(\bar\phi(t),\bar\pp(t))$
which are exponentially decreasing as well, with the same exponent $\alpha$.

Recall also that the Lipshitz constant of $\qq$ is uniformly small (and tends
to zero as $\eb\to0$ and $\gamma\to0$). The Lipshitz constant of $\tilde{\Bbb V}$ with respect to $\varphi$
is uniformly bounded and the Lipshitz constant
with respect to $\pp$ is uniformly small (and tends to zero as $\eb\to0$;
see \eqref{tlvtlu},\eqref{wman},\eqref{2.25}). The Lipshitz constant of $\hh$ with respect to $\pp$ is of
order $\eb$ and the Lipshitz constants with respect to $\varphi$ and $\vv$
are bounded (see\eqref{hkforno},\eqref{fgcr}). Thus, for the functions under the integrals in
\eqref{2.cor}, the Lipshitz
constants with respect to both $\varphi$ and $\pp$ in the first equation
of \eqref{2.cor} and with respect to $\pp$ in the second equation
are uniformly small, while the Lipshitz constant with respect to $\varphi$
in the second equation is uniformly bounded. This immediately implies that
operator \eqref{2.cor} is contracting on the space of exponentially decreasing
functions $(\phi(t),\pp(t))_{t\geq0}$ endowed with the norm
$\displaystyle \|\phi,\pp\|=\sup_{t\geq0} e^{\alpha t} (\|\varphi(t)\|_\Psi+\kappa\|\pp(t)\|_{\Bbb P})$,
for all sufficiently small $\kappa$, $\eb$ and $\gamma$. This gives us the
required existence and uniqueness of the fixed point $(\phi(t),\pp(t))_{t\geq0}$.
Being the fixed point of a contracting operator, it depends continuously on
every parameter on which the operator depends continuously, so it
depends continuously on the function $\tilde{\Bbb V}^2$. Note that
$\tilde{\Bbb V^2}\to\tilde{\Bbb V^1}$ as $\gamma\to0$ (by \eqref{tlvtlu},
this just means that the manifold $\Bbb W_{\yy^0,\eb}$
depends on $\yy^0$ continuously). Hence, in the same limit
we have $(\phi(t),\pp(t))\to 0$ (which is the trivial fixed point of \eqref{2.cor}
when $\tilde{\Bbb V^2}\equiv\tilde{\Bbb V^1}$).
\end{proof}

\section{Spatially non-walking solutions and their entropy}\label{s3}

In our application to Ginzburg-Landau equation,
the $\pp$-component of the LDS \eqref{2.10} describes the
temporal evolution of the centers of soliton pairs, namely the deviations of the pair centers
from the points of a given spatial lattice. This description
is valid only if the distances between the soliton pairs are large enough, i.e. the deviations
of the soliton pairs from the lattice points stay uniformly bounded for arbitrarily large
lattice sizes, see Section \ref{s4}.
Thus, it is crucial to be able to control the norm of $\pp(t)=\{p_k(t)\}_{k=-\infty}^{+\infty}$,
i.e. to keep all $p_k$ bounded.

On the other hand, according to \eqref{pp1e},\eqref{ph1e},
in the zero order approximation with respect to $\eb$ we have
\begin{equation}\label{3.1}
p_k(t)\approx p_k(0)+\int_0^tg(y^0_k(s))\,ds,
\end{equation}
where $y^0_k$ is a trajectory from the given hyperbolic set $\Lambda$.
Thus, an independent diffusive-like behavior of the coordinates
$p_k(t)$ should be expected \cite{46} in the case $\Lambda$ is non-trivial (chaotic),
i.e. the quantities $p_k(t)$ are out of control in this case.

The main aim of the Section is to show, however, that under some natural
assumptions on the set $\Lambda$ there exists a
set of solutions for which at all $t\in\Bbb R$
\begin{equation}\label{3.2}
\|\pp(t)\|_{_{\Bbb P}}\le R_0
\end{equation}
for some constant $R_0>>1$. Moreover, this set is large enough, so that it has
positive space-time entropy. In what follows, in order to simplify notations,
we assume that all individual ODE's in the uncoupled
LDS \eqref{2.7} are identical, i.e. $f_k\equiv f$, $g_k\equiv g$ for all $k$.

\begin{theorem}\label{Th3.1} Let the assumptions of Theorem \ref{Th2.1} hold
and let $N:=\dim P$. Let us also assume that the hyperbolic set $\Lambda$ of system
\eqref{2.1} is transitive and locally-maximal and contains $N+1$ periodic orbits $Z_1: y=z_1(t)$,
$Z_2: y=z_2(t)$,$\dots$, $Z_{N+1}: y=z_{N+1}(t)$ with periods $T_1$,$\dots$, $T_{N+1}$
respectively. Define the vectors $\vec b_i\in P$, $i=1,\dots,N+1$, as follows:
\begin{equation}\label{3.21}
 \vec b_i:=\frac1{T_i}\int_0^{T_i}g(z_i(t))\,dt,
\end{equation}
and require the following properties to be satisfied:

1. linear combinations of vectors $\vec b_i$ generate the whole space $P$:
\begin{equation}\label{span}
P=\sspan\{\vec b_1,\dots,\vec b_{N+1}\};
\end{equation}

2. there exist strictly positive numbers $A_i$ such that
\begin{equation}\label{3.22}
A_1\vec b_1+A_2\vec b_2+\dots+A_{N+1}\vec b_{N+1}=0.
\end{equation}
Then, for all sufficiently small $\eb>0$, there exists a uniformly bounded set $\Cal K$ of solutions
of system \eqref{2.10} which has strictly positive space-time entropy:
\begin{equation}\label{3.33}
h(\Cal K)>0.
\end{equation}
\end{theorem}
\begin{proof} We start with the following observation.
\begin{lemma}\label{Lem3.1} Let (\ref{span}) and (\ref{3.22}) hold. Then, for every vector $p\in P$, $p\ne0$,
there exists $j=J(p)\in\{1,\cdots, N+1\}$ such that
\begin{equation}\label{3.34}
p \cdot \vec b_{J(p)}<0
\end{equation}
and, consequently,  there exists $\delta>0$ such that, for every $p\ne0$,
\begin{equation}\label{3.35}
\cos(p,\vec b_{J(p)})\le-\delta.
\end{equation}
\end{lemma}
\noindent Indeed, suppose there exists $p$ such that \eqref{3.34} is wrong, i.e. $p\cdot \vec b_i\ge0$ for all
$i=1,\dots, N+1$. Multiplying then equality \eqref{3.22}
by this $p$ and using that $A_i>0$, we conclude that $p\cdot\vec b_i=0$ for all $i$. By \eqref{span},
this contradicts the assumption $p\neq 0$.
Thus, \eqref{3.34} is verified and \eqref{3.35} follows immediately from \eqref{3.34} by compactness arguments.

The idea of the proof of the theorem is as follows. As the set $\Lambda$ is transitive and locally-maximal,
for every two of the periodic orbits $Z_i$ and $Z_j$ we may choose two different heteroclinic orbits
$Z_{ijm}: y=z_{ijm}(t)$, $m=1,2$, that connect them, i.e.
$$\lim_{t\rightarrow-\infty} (z_{ijm}(t)-z_i(t+\theta^-_{ijm}))=0,\qquad
\lim_{t\rightarrow+\infty} (z_{ijm}(t)-z_j(t+\theta^+_{ijm}))=0$$
for some constant $\theta^{\pm}_{ijm}$. The orbits $Z_{ijm}$ also belong to $\Lambda$; in fact, the number of
different heteroclinics is infinite for each pair of periodic orbits in $\Lambda$, but we
need only two of them for each $i$ and $j$. The existence of the heteroclinics mean that we may build
orbits in $\Lambda$ which stay for some time near the orbit $Z_i$, then "jump" along any two of the heteroclinics
$Z_{ij1,2}$ into a neighborhood of $Z_j$, stay there, then jump again into a neighborhood of another
periodic orbit, etc.. We will see that for sufficiently small $\eb$ one can build orbits
$\yy=\{y_k(t)\}_{k=-\infty}^{k=+\infty}$ of system \eqref{2.10} with a similar behavior for every
component $y_k(t)$: the component stays close to $z_i(t)$ for some time then jumps to $z_j(t)$, etc.,
moreover the choice of the sequence of the periodic orbits the component shadows can be
made independently for different $k$. When the component $y_k$ is close to $z_i(t)$ for sufficiently
long time, the $p_k$-component of the associated solution will move in the direction close to $\vec b_i$
as time grows (see \eqref{3.1},\eqref{3.21}). By \eqref{3.35}, if the norm of $p_k$ becomes large enough
we can always find a vector $\vec b_j$ such that moving in its direction will lead to a decrease in the norm
of $p_k$. Thus, by jumping each time to a properly chosen periodic orbit $Z_j$ we may keep the norm of
all $p_k$ bounded. As each jump can be made by at least two different ways (along the first or the second
heteroclinic) the set of different solutions of system \eqref{2.10} we obtain in this way will have positive
entropy.

As the first step in implementing this construction we recall the following
standard result on the ``shadowing'' in hyperbolic sets.

\begin{lemma}\label{Lem3.etr} There exist $\bar\gamma>0$ and $\alpha>0$
such that for any two orbits $y_-(t)$ and $y_+(t)$ from the hyperbolic
set $\Lambda$ which satisfy
\begin{equation}\label{3.close}
\|y_-(t_0)- y_+(t_0)\|\le\gamma,
\end{equation}
where $\bar\gamma>\gamma>0$, there exists an orbit $y(t)\in\Lambda$ and a phase shift $\theta$ such that
\begin{equation}\label{3.tr}
\begin{array}{l}
\|y(t)-y_-(t)\|\le C_\gamma e^{\alpha t} \;\mbox{ for } t\le t_0, \\ \\
\|y(t)-y_+(t+\theta)\|\le C_\gamma e^{-\alpha t} \;\mbox{ for } t\ge t_0,\qquad
|\theta|\le C_\gamma,
\end{array}
\end{equation}
where $C_\gamma>0$ depends only on $\gamma$ and tends to zero as $\gamma\to 0$.
\end{lemma}
A proof can be found e.g. in \cite{23}. The orbit $y(t)$ corresponds simply to the intersection
of the local unstable manifold of $y_-(t)$ with the local stable manifold of $y_+(t)$; this intersection
belongs to $\Lambda$ because this set is locally-maximal.

Combining Lemma \ref{Lem3.etr} with Theorem \ref{Prop2.trace}, we obtain an analogous result
for the lattice dynamical system \eqref{2.10}.
\begin{lemma} \label{Lem3.fin} There exist $\alpha>0$ and $\bar\gamma>0$ such that for all $\eb>0$ small enough,
for any two orbits $\yy^0_\pm\in\Lambda^\infty$ such that
\begin{equation}\label{3.cl}
\|\yy^0_-(t_0)-\yy^0_+(t_0)\|_{_{\Bbb Y}}\le\gamma,
\end{equation}
where $\bar\gamma>\gamma>0$ and $t_0\in\R$, and for any solution $(\yy_-(t), \pp_-(t))$ of \eqref{2.10}
belonging to the invariant manifold ${\Bbb W}_{\yy^0_-,\eb}$, there exist an
orbit $\yy^0\in\Lambda^\infty$, a solution
$(\yy(t),\pp(t))\in {\Bbb W}_{\yy^0,\eb}$ of the lattice system \eqref{2.10}, and the set
of constant phase shifts $\theta_k$, $k\in\Bbb Z$, such that
\begin{equation}\label{3.np00}
\begin{array}{l}
\sup_{k\in{\Bbb Z}} \|y^0_k(t)-y^0_{k+}(t+\theta_k)\|_Y\le C_\gamma
e^{-\alpha (t-t_0)},\; \mbox{ for } t\ge t_0,\\ \\
\sup_{k\in{\Bbb Z}} \|y^0_k(t)-y^0_{k-}(t)\|_Y\le C_\gamma
e^{\alpha(t-t_0)} \; \mbox{ for } t\le t_0,\qquad\;\;
\sup_{k\in{\Bbb Z}}\|\theta_k\|\le C_\gamma
\end{array}
\end{equation}
and
\begin{equation}\label{3.neg}
\|\yy(t)-\yy_-(t)\|_{_{\Bbb Y}}+\|\pp(t)-\pp_-(t)\|_{_{\Bbb P}}\le
 C_\gamma\;e^{\alpha(t-t_0)}, \ \ t\le t_0,\\
\end{equation}
where $C_\gamma\to+0$ as $\gamma\to 0$.
\end{lemma}

\begin{proof} Indeed, in order to find the required solution $(\yy(t),\pp(t))$, we first
construct a trajectory $\yy^0\in\Lambda^\infty$, each component $y_k(t)$ of which
is defined by $y^0_{k-}(t)$ and $y^0_{k+}(t)$ by virtue of Lemma \ref{Lem3.etr} such that
\eqref{3.np00} is satisfied (since the
unperturbed system \eqref{2.7} is a Cartesian product of systems
\eqref{2.1}, we only need to apply Lemma \ref{Lem3.etr} to every component
in this product). Applying after that Remark
\ref{plm}, we find (in a unique way) the solution $(\yy(t),\pp(t))$ of the
perturbed system \eqref{2.10}, satisfying \eqref{3.neg} for $t\le t_0$.
\end{proof}

We are now ready to complete the proof of the theorem. We will choose sufficiently large constants $T$ and
$R$ and sufficiently small constants $\nu$ and $\mu$ and construct a sequence of
sets ${\Cal K}_l$ of solutions of \eqref{2.10} and a sequence
of sets ${\Cal K}^0_l$ of orbits from $\Lambda^\infty$ such that:\\
1) for each of the solutions from ${\Cal K}_l$ there exists an otbit $\yy^0\in{\Cal K}^0_l$
such that the solution belongs to the invariant manifold ${\Bbb W}_{\yy^0,\eb}$;\\
2) for every trajectory $\yy^0=\{y_k(t)\}_{k=-\infty}^{+\infty}\in {\Cal K}^0_l$,
for every $k\in {\Bbb Z}$ there are periodic orbits $Z_{i_{k+}}: y=z_{i_{k+}}(t)$ and
$Z_{i_{k-}}: y=z_{i_{k-}}(t)$ (from the set of
periodic orbits $Z_1,\dots,Z_{N+1}$ under consideration) such that
\begin{equation}\label{cond2}
\|y_k^0(lT)-z_{i_{k+}}(\tau_{kl})\|< \nu, \qquad \|y_k^0(-lT)-z_{i_{k-}}(\tau_{k,-l})\|< \nu
\end{equation}
for some (irrelevant) constants $\tau\in[0,{\Cal T}]$, where $\Cal T=\max_{i=1,\dots,N+1} T_i$
(the periods of $Z_i$);\\
3) at $t_0=\pm lT$
\begin{equation}\label{3.ind}
\|\pp(t_0)\|_{_{\Bbb P}}\le R;
\end{equation}
4) for every solution $(\tilde\yy(t),\tilde\pp(t))\in {\Cal K}_{l+1}$
there exists a solution $(\yy(t),\pp(t))\in {\Cal K}_l$ such that
\begin{equation}\label{cond4}
\|\tilde\yy(t)-\yy(t)\|_{_{\Bbb Y}}+\|\tilde\pp(t)-\pp(t)\|_{_{\Bbb P}}\le
\mu e^{-\alpha(lT-|t|)}
\end{equation}
for all $|t|\leq lT$ (the constant $\alpha>0$ depends on the set $\Lambda$ only).

By condition 4, the sequence of the sets ${\Cal K}_l$ converges, as $l\rightarrow+\infty$,
to a certain set $\Cal K$ of solutions of the LDS \eqref{2.10} (convergence is uniform on
any bounded time interval). Moreover, conditions 2, 3 and 4 imply that solutions in
the set $\Cal K$ are uniformly bounded, in particular $\|\pp(t)\|$ is uniformly bounded for all of
the solutions. Thus, to prove the theorem, we need to actually construct the sequence
${\Cal K}_l$ and to do it in such a way that the sets ${\Cal K}_l$ would contain
``sufficiently many'' solutions --
this would warrant the positivity of the space-time entropy of the limit set $\Cal K$.

As ${\Cal K}^0_0$ we choose the set that consists of one orbit
$\yy^0(t)=\{y_k^0=z_1(t)\}_{k=-\infty}^{\infty}$; the set ${\Cal K}_0$ will
consist of one solution in the invariant manifold ${\Bbb W}_{\yy^0,\eb}$ which satisfies
$\pp(0)=0$.

Now assume we have built the sets ${\Cal K}^0_l$, ${\Cal K}_l$ for some $l$, and let us construct
the sets ${\Cal K}^0_{l+1}$, ${\Cal K}_{l+1}$. Take any pair
$\left\{\yy^0\in{\Cal K}^0_l,\;(\yy,\pp)\in{\Cal K}_l\cap{\Bbb W}_{\yy^0,\eb}\right\}.$
Let $i_{k\pm}$ ($k\in{\Bbb Z}$)
be the sequences of indices defined by \eqref{cond2} and $j_{k\pm}:=J(p_k(\pm lT))$, where the integer-valued
function $J(p)$ is defined by \eqref{3.35}. Choose any two sequences $m_{k\pm}$ ($m_{k\pm}=1$ or $2$).
Choose an
orbit $\yy^0_+\in\Lambda^\infty$ as follows: $\yy^0_{k+}(t)=z_{i_{_{k+}}j_{_{k+}}m_{_{k+}}}(t-lT+\tau_{kl})$,
where $y=z_{ijm}(t)$
is one of the two (chosen above) heteroclinic orbits $Z_{ij1,2}$ which connect the periodic orbits
$Z_i$ and $Z_j$. We assume here that the time parametrization on the heteroclinic orbits is chosen such that
$\|z_{ijm}(t)-z_i(t)\|=\nu$ at $t= {\Cal T}$, and $\|z_{ijm}(t)-z_i(t)\|<\nu$ at all $t< {\Cal T}$. Hence,
$\|z_{ijm}(\tau_{kl})-z_i(\tau_{kl})\|\leq \nu$ (recall that the numbers
$\tau_{kl}$ are bounded by $\Cal T$), so
$\|\yy^0-\yy^0_+\|_{_{\Bbb Y}}< 2\nu$ by \eqref{cond2}. Therefore, if $\nu$ is small enough, we may apply
Lemma \ref{Lem3.etr} (with $\yy^0$ taken as the orbit $\yy^0_-$ of the lemma) and obtain a solution
$(\hat \yy(t), \hat \pp(t))$ such that
\begin{equation}\label{cond4hat}
\|\hat\yy(t)-\yy(t)\|_{_{\Bbb Y}}+\|\hat\pp(t)-\pp(t)\|_{_{\Bbb P}}\leq \mu e^{-\alpha(lT-t)}
\end{equation}
at $t\leq t_0=lT$; moreover this solution belongs to the invariant manifold ${\Cal W}_{\hat\yy^0,\eb}$
associated with the orbit $\hat\yy^0\in\Lambda^\infty$ such that, as $t\rightarrow+\infty$, the components
$\hat y_k^0(t)$ tend exponentially
to the heteroclinic orbits $Z_{i_{_{k+}}j_{_{k+}}m_{_{k+}}}$ -- hence to the periodic orbits $Z_{j_{k+}}$.

Absolutely analogously (by applying the version of Lemma \ref{Lem3.etr} obtained by inversion of time) we
obtain the existence of a solution $(\tilde \yy(t), \tilde \pp(t))$ such that
\begin{equation}\label{cond4tild}
\|\tilde\yy(t)-\hat\yy(t)\|_{_{\Bbb Y}}+\|\tilde\pp(t)-\hat\pp(t)\|_{_{\Bbb P}}\leq \mu e^{-\alpha(t+lT)}
\end{equation}
at $t\geq t_0=-lT$; moreover this solution belongs to the manifold ${\Cal W}_{\tilde\yy^0,\eb}$
associated with the orbit $\tilde\yy^0\in\Lambda^\infty$ such that at each $k$ the component
$\tilde y_k^0(t)$ tends exponentially to the heteroclinic orbit
$Z_{j_{_{k-}}i_{_{k-}}m_{_{k-}}}$ as $t\rightarrow-\infty$ (and it still tends to
$Z_{i_{_{k+}}j_{_{k+}}m_{_{k+}}}$
as $t\rightarrow+\infty$).

By \eqref{cond4tild}, \eqref{cond4hat}, condition \eqref{cond4} is fulfilled by the newly constructed
solution $(\tilde\yy,\tilde\pp)$. Since each component $\tilde y^0_k(t)$ tends to the corresponding
periodic orbit $Z_{j_{k+}}$ as $t\rightarrow+\infty$ and to $Z_{j_{k-}}$ as $t\rightarrow-\infty$,
and the convergence is, by construction, uniform for all $k$, $l$ and for all possible initial solutions
$(\yy(t),\pp(t))\in{\Cal K}_l$, it follows that condition \eqref{cond2}
will be satisfied for the orbit $\tilde \yy^0$ at $t=\pm (l+1)T$, provided $T$ was chosen large enough.

It follows, furthermore, that if $T$ is sufficiently large and $\eb$ is sufficiently
small, then the change in $\tilde p_k$ along the orbit $(\tilde \yy(t),\tilde \pp(t))$ for the time from
$t=lT$ to $t=(l+1)T$ equals to $T\vec b^\prime_{j_{k+}}$ where $\vec b^\prime_{j_{k+}}$ is uniformly close to
the vector $\vec b_{j_{k+}}$ defined by \eqref{3.21}. As $j_{k+}=J(p_k(lT))$ and $\tilde p_k(lT)$ is
close to $p_k(lT)$ (see \eqref{cond4}), it follows that
$$
\cos(\tilde p_k(lT),\vec b^\prime_{j_{k+}})<-\delta/2, \ \ k\in\Bbb Z
$$
(see \eqref{3.35}). Therefore,
\begin{equation}\label{3.back}
\tilde p_k^2((l+1)T)=(\tilde p_k(lT)+T\vec b^\prime_{j_{k+}})^2\leq \tilde p_k^2(lT) +
T \|\vec b^\prime_{j_{k+}}\| ( T \|\vec b^\prime_{j_{k+}}\|-\delta \|\tilde p_k\|).
\end{equation}
By \eqref{cond4} and \eqref{3.ind}
$\|\tilde p_k(lT)\|\leq R+\mu$, and we see now from \eqref{3.back} that
$$\tilde p_k^2((l+1)T)<R^2,$$
provided $T$ is taken sufficiently large with respect to $\mu$ and $R$ is sufficiently large with respect to $T$
(note that $\|\vec b_j\|$ is bounded away from zero by virtue of \eqref{span}, \eqref{3.22}).
Analogously, one checks that
$$\tilde p_k^2(-(l+1)T)<R^2.$$
As we see, condition \eqref{3.ind} is satisfied by the solution $(\tilde\yy,\tilde\pp)$ at $t_0=\pm (l+1)T$.

Thus, we have shown that given any solution from the set ${\Cal K}_{l_0}$ and any pair of sequences
$m_{k\pm}$ (these sequences define which of the two heteroclinic connections is used to jump from
the periodic orbit $Z_{i_{k\pm}}$ to $Z_{j_{k\pm}}$) we obtain a solution which satisfies above
conditions 1-3 with $l=l_0+1$, i.e. the newly built solution can be included into the set ${\Cal K}_{l_0+1}$;
we have also checked condition 4 that ensures the convergence of the sequence of sets ${\Cal K}_l$
as $l\rightarrow+\infty$. As we may choose the sequences $m_{k\pm}$ in an arbitrary way at each
step of the procedure, the number of solutions in the set ${\Cal K}_l$ which stay at a bounded away from zero
distance from each other at $|t|\leq lT$ and $|k|\leq n$ equals to $4^{l(2n+1)}$. This immediately shows
that the space-time entropy of the limit set $\Cal K$ is strictly positive.
\end{proof}

Note that the assumption that the set $\Lambda$ is locally-maximal and transitive
can be formulated in a more constructive way. Indeed, assume that we have a set
of hyperbolic periodic orbits $Z_1,\dots,Z_{N+1}$, which
satisfy conditions 1 and 2 of the theorem. Build an oriented graph with $N+1$ vertexes:
the edge connects the vertex $i$ with vertex $j$ if we know there exists a heteroclinic orbit
$Z_{ij}$ which corresponds to a {\em transverse} intersection of the unstable manifold of $Z_i$ with
the stable manifold of $Z_j$. If this graph is transitive, then the set $\Lambda$
of all orbits which stay for all times in a sufficiently small neighborhood of
the union of the hyperbolic periodic orbits $Z_i$ and the transverse heteroclinic
orbits $Z_{ij}$ is uniformly-hyperbolic, transitive and locally-maximal
\cite{AS}, so Theorem \ref{Th3.1} holds.

Note also that assumption \eqref{3.22} is really important for the
proof of the theorem. Indeed, consider the case $\dim P=1$ for example. Here the
integrals $\vec b_1$ and $\vec b_2$ are real numbers, and if condition \eqref{3.22}
is violated, they both have the same sign, positive, say.
In this case, when the component $y_k$ stays close to either of the
periodic orbits $Z_{1,2}$, the component $p_k$ will increase with time, so we cannot
keep $p_k(t)$ bounded by mere switching between $Z_1$ and $Z_2$. However,
assumption \eqref{3.22} can be relaxed if we allow for a uniform drift, common for
all $p_k(t)$. Namely, the following statement holds true.

\begin{corollary}\label{Cor3.drift} Let all of the assumptions of
Theorem \ref{Th3.1} be fulfilled except of \eqref{3.22}. Assume
the convex hull $I_b$ of vectors $\vec b_1,\, \cdots\, \vec b_{N+1}$ have a non-empty
interior:
\begin{equation}\label{3.pb}
I_b:=\operatorname{int}\{\operatorname{conv}\{\vec b_1,\cdots,\vec b_{N+1}\}\}\ne\varnothing.
\end{equation}
Let $\vec p\subset I_b$. Then, for every sufficiently small
$\eb$, there exists a set $\Cal K_{\vec p}$ of solutions $(\yy(t),\pp(t))$ of
system \eqref{2.10} such that $\Cal K_{\vec p}$ has positive space-time entropy and each solution from
$\Cal K_{\vec p}$ satisfies
\begin{equation}\label{3.drboun}
\|\yy(t)\|_{_{\Bbb Y}}+\|\pp(t)-\vec p t\|_{_{\Bbb P}}\le R_0<\infty,\ \ t\in\R,
\end{equation}
where the constant $R_0$ depends on $\vec p$, but is independent of $t$ and the choice of
the solution.
\end{corollary}
\noindent Indeed, for every $\vec p\in I_b$  conditions \eqref{span},\eqref{3.22} hold
for the vectors $\vec b_1-\vec p$, $\vec b_2-\vec p$, $\dots$, $\vec b_{N+1}-\vec p$.
Then, applying Theorem \ref{Th3.1} to the system obtained from \eqref{2.10} by substracting
$\vec p$ from the function $g$, we immediately obtain the corollary.

\begin{remark}\label{Rem3.1d}
In the one-dimensional case ($\dim P=1$), we only
need two hyperbolic periodic orbits, $Z_+$ and $Z_-$, connected by transverse
heteroclinics. Conditions \eqref{span}, \eqref{3.22} read now
\begin{equation}\label{3.1dll}
\int_0^{T_-}g(z_-(t))\,dt \; \cdot \; \int_0^{T_+}g(z_+(t))\,dt <0;
\end{equation}
conditions \eqref{span}, \eqref{3.pb} read as
\begin{equation}\label{3.1DOK}
\frac{1}{T_-}\int_0^{T_-}g(z_-(t))\,dt\ne \frac{1}{T_+}\int_0^{T_+}g(z_+(t))\,dt.
\end{equation}
In order to establish the existence of the heteroclinic cycle with two hyperbolic
periodic orbits one may use Shilnikov criterion. Namely, it is enough to show the
existence of a saddle-focus equilibrium state $y=z_0$ with a homoclinic loop $y=z_h(t)$,
$z_h(t)\rightarrow z_0$ as $t\rightarrow\pm\infty$,
and to check that the so-called
Shilnikov conditions of chaos are satisfied (we will not discuss a higher-dimensional case
as in the application we consider in this paper we have $y\in R^3$; in the three-dimensional
case the Shilnikov condition is that the nearest to the imaginary axis characteristic
exponent is not real; the equilibrium state must be hyperbolic, i.e. it has
characteristic exponents on both sides of the imaginary axis and no characteristic
exponents on the axis). Then there exists a sequence $Z_n$ of hyperbolic periodic orbits
which converge to the homoclinic loop as $n\rightarrow+\infty$, any two of them are
connected by transverse heteroclinics \cite{Sh1,Sh2}. The periods $T_m$ of $Z_m$ tend to infinity.
One can always choose time parametrization such that
$\sup_{t\in[-\frac{T_m}{2},\frac{T_m}{2}]} |z_m(t)-z_h(t)-z_0|\rightarrow 0$
as $m\rightarrow+\infty$. It follows that one can always choose among the orbits $Z_m$
a pair satisfying conition \eqref{3.1DOK}, provided
\begin{equation}\label{sfcond}
\int_{-\infty}^{+\infty}(g(z_h(t))-g(z_0))\,dt\ne 0
\end{equation}
(the integral converges since $z_h(t)$ tends to $z_0$ exponentially - because of
the hyperbolicity of $z_0$). Note that the homoclinic loop to a saddle-focus may split
as we perturb the system, however the two hyperbolic
periodic orbits that we find near the loop do not disappear, nor
the transverse heteroclinics that connect them do, so by checking condition
\eqref{sfcond} for one parameter value we establish the existence of spatio-temporal chaos
for an open set of parameter values; see Lemma \ref{Lemzero3} in Section \ref{s4} for an example.
\end{remark}

\begin{remark}\label{Rem3.MD} In the case of LDS with $n$ spatial dimensions, i.e.
those parameterized by multiindices
$k\in\Bbb Z^n$ instead of $k\in\Bbb Z$,
the result of Theorem \ref{Th3.1},
obviously, remains true under the properly modified definition
of the space-time topological entropy. In fact, this case is just formally
reduced to $k\in \Bbb Z$ by an appropriate reparameterization of the grid $\Bbb Z^n$
by the points from $\Bbb Z$.
\end{remark}
\section*{Appendix. Soliton interaction equations}
{\tiny We do not reproduce a proof here, however
we provide a very informal derivation of the key equations \eqref{4.lrde}, for the sake of completeness. Write
equation \eqref{4.ppcGL} in a symbolic form $\Dt u={\Cal H}(u)+\eb \Cal G(u)$,
where the operator $\Cal H$ is such that $\Cal H(U)=0$, and the small parameters $\eb$
govern the perturbation, i.e. $\eb=(\mu,\delta-\delta_0)$ in our case, and
$\eb \Cal G(u):=\mu-i(\delta-\delta_0)u$. By making the multi-pulse ansatz
$u\approx\sum U_{\xi_j(t),\phi_j(t)}$ (see \eqref{4.mp}), we obtain the following equation:
$\sum \{-e^{i\phi_j}\varphi_1(x-\xi_j)\dot\xi_j+e^{i\phi_j}\varphi_2(x-\xi_j)\dot\phi_j\}
\approx{\Cal H}(u)+\eb \Cal G(u)$ (see \eqref{4.man1}),
which gives us, after taking the inner product (defined by \eqref{normaf})
with $e^{-i\phi_j}\psi_{1,2}(x-\xi_j)$, that
$\dot \xi_j\approx -\Ree \int_{-\infty}^{+\infty}e^{-i\phi_j}\psi_{1}(x-\xi_j)\{\Cal H+\eb \Cal G\}dx$, and
$\dot \phi_j\approx \Ree \int_{-\infty}^{+\infty}e^{-i\phi_j}\psi_{2}(x-\xi_j)\{\Cal H+\eb \Cal G\}dx$.
As the functions $\psi_{1,2}$ decay exponentially, and the distances $|\xi_i-\xi_{i-1}|$ assumed all to be
large, the main contribution to these integrals is given by the segment
$x\in[\frac{1}{2}(\xi_{j-1}+\xi_j),\frac{1}{2}(\xi_{j}+\xi_{j+1})]$. As the pulse $U(x)$ decays exponentially
fast, the main contribution to the approximate solution
$u=\sum U_{\xi_k,\phi_k}$ on this segment is given by the term $U_{\xi_j,\phi_j}:=e^{i\phi_j}U(x-\xi_j)$,
and the next order corrections come from the two terms $U_{\xi_{j\pm1},\phi_{j\pm1}}$. By recalling that
$\Cal H(U_{\xi_j,\phi_j})=0$, we obtain (after shifting $x\mapsto x+\xi_j$ in the integral)
$\dot \xi_j\approx -\Ree \int_{-\infty}^{+\infty}\theta_j(x)\psi_{1}(x)\Cal L_{_U}
\{U_{\xi_{j+1}-\xi_j,\phi_{j+1}-\phi_j}+U_{\xi_{j-1}-\xi_j,\phi_{j-1}-\phi_j}\}dx
+\eb\Ree \int_{-\infty}^{+\infty}e^{-i\phi_j}\psi_{1}(x) \Cal G(e^{i\phi_j}U(x))dx$, where
$\Cal L_{U}:=\Cal H'(U)$ is given by \eqref{lulin}, and $\theta_j(x)=1$
at $x\in[\frac{1}{2}(\xi_{j-1}-\xi_j),\frac{1}{2}(\xi_{j+1}-\xi_{j})]$ and vanishes at the other
values of $x$. We may rewrite the equation for $\dot\xi_j$ as
$\dot \xi_j\approx -\Ree \int_{-\infty}^{+\infty}\Cal L_{_U}^\dagger\{\theta_j(x)\psi_{1}(x)\}
(U_{\xi_{j+1}-\xi_j,\phi_{j+1}-\phi_j}+U_{\xi_{j-1}-\xi_j,\phi_{j-1}-\phi_j})dx
+\eb\Ree \int_{-\infty}^{+\infty}e^{-i\phi_j}\psi_{1}(x) \Cal G(e^{i\phi_j}U(x))dx$, where
$\Cal L_{_U}^\dagger$ is conjugate to $\Cal L_{_U}$ (see \eqref{lulc}). As $\Cal L_{_U}^\dagger \psi_{1}=0$,
it is easy to see from \eqref{lulc} that $\Cal L_{_U}^\dagger\{\theta_j(x)\psi_{1}(x)\}=
(1+i\beta_0)[\theta_j''(x)\psi_{1}(x)+2\theta_j'(x)\psi_{1}'(x)]$. By noticing that $\theta_j'(x)$
is the difference of two delta-functions, the one at
$x=\kappa_j^-:=\frac{1}{2}(\xi_{j-1}-\xi_j)$ minus the other at $x=\kappa_j^+:=\frac{1}{2}(\xi_{j+1}-\xi_{j})$,
one immediately
evaluates the first integral in the equation for $\dot\xi$; to find that the second integral is zero
we notice that
in our case $\eb e^{-i\phi_j}\Cal G(e^{i\phi_j}U(x))\equiv \mu e^{-i\phi_j}-(\delta-\delta_0)\varphi_2$
and use the normalization condition \eqref{normaf}. As the result we get\\
$\dot \xi_j\!\approx\Ree(1\!+ i\beta_0)\!\left\{\!e^{i(\phi_{j-1}\!-\phi_j)}[\psi_{1}(x)
U'(x\!-\!2\kappa_j^-)\!-\psi_{1}'(x)U(x\!-\!2\kappa_j^-)]_{x=\kappa_j^-}-
e^{i(\phi_{j+1}\!-\phi_j)}[\psi_{1}(x)U'(x-\!2\kappa_j^+)\!-
\psi_{1}'(x)U(x-\!2\kappa_j^+)]_{x=\kappa_j^+}\!\!\right\}$.
Similarly,\\
$\dot \phi_j\!\approx-\Ree(1\!+ i\beta_0)\!\left\{\!e^{i(\phi_{j-1}\!-\phi_j)}[\psi_{2}(x)
U'(x\!-\!2\kappa_j^-)\!-\psi_{2}'(x)U(x\!-\!2\kappa_j^-)]_{x=\kappa_j^-}-
e^{i(\phi_{j+1}\!-\phi_j)}[\psi_{2}(x)U'(x-\!2\kappa_j^+)\!-
\psi_{2}'(x)U(x-\!2\kappa_j^+)]_{x=\kappa_j^+}\!\!\right\}$

$+\mu\Ree ce^{i(\zeta-\phi_j)} -(\delta-\delta_0)$\\
(see \eqref{normaf},\eqref{4.cc}). As $\kappa_j^\pm=\frac{1}{2}(\xi_{j\pm1}-\xi_j)$ are large in absolute value,
we may use asymptotics \eqref{4.tails}, which gives us expressions for $\dot\xi_j$ and $\dot\phi_j$
compatible with \eqref{4.rde} indeed. Similar derivations can be found e.g.
in \cite{46}; a complete proof is in \cite{57}.}

\end{document}